\documentclass[10pt]{article}
\usepackage{times,amsmath,amsthm,amssymb,graphicx,xspace,epsfig,xcolor,colortbl,diagbox}

\usepackage{algorithm}\usepackage{algorithmic}
\usepackage{tikz, tkz-graph, tkz-berge}
\usetikzlibrary{decorations.pathreplacing}
\usetikzlibrary{patterns}
\usetikzlibrary{patterns.meta}
\usepackage{pgf}
\usepackage{pgffor}
\usepackage{xspace}
\usepackage{color}
\usepackage{comment}
\usepackage{authblk}
\usepackage{multirow}
\usepackage{hyperref}
\usepackage{cleveref}
\usepackage{dsfont}
\usepackage{mathtools,amsfonts,amsthm,mathrsfs,amssymb,stmaryrd}
\usepackage{enumitem}
\usepackage{thm-restate}

\usepackage{caption}
\usepackage{subcaption}
\usepackage{soul,cite}

\newtheorem{theorem}{Theorem}
\newtheorem{corollary}[theorem]{Corollary}
\newtheorem{proposition}[theorem]{Proposition}
\newtheorem{lemma}[theorem]{Lemma}

\theoremstyle{definition}
\newtheorem{remark}[theorem]{Remark}
\newtheorem{question}[theorem]{Question}

\newtheorem{conjecture}[theorem]{Conjecture}

\newtheorem{claim}[theorem]{Claim}

\newenvironment{proofclaim}[1][]
	{\par\noindent {\it Proof of the claim}. }{ \hfill$\lozenge$\par\vspace{11pt}}

\DeclareMathOperator{\Mad}{mad}
\DeclareMathOperator{\Ad}{ad}
\DeclareMathOperator{\UG}{UG}
\DeclareMathOperator{\diam}{diam}
\DeclareMathOperator{\tw}{tw}
\DeclareMathOperator{\dtw}{dtw}

\DeclareMathOperator{\girth}{girth}
\DeclareMathOperator{\dw}{\mathscr{D}w}

\definecolor{lightgrey}{RGB}{225, 225, 225}

\newcommand{\cev}[1]{\reflectbox{\ensuremath{\vec{\reflectbox{\ensuremath{#1}}}}}}
\newcommand{\dic}{\vec{\chi}}
\newcommand{\bid}{\overleftrightarrow}

\begin{document}

\title{Redicolouring digraphs: directed treewidth and cycle-degeneracy\thanks{Research supported by the CAPES-Cofecub project Ma 1004/23, by the Inria Associated Team CANOE, by the research grant
    DIGRAPHS ANR-19-CE48-0013 and by the French government, through the EUR DS4H Investments in the Future project managed by the National Research Agency (ANR) with the reference number ANR-17-EURE-0004.}}
\author[1]{Nicolas Nisse}
\author[1]{Lucas Picasarri-Arrieta}
\author[2]{Ignasi Sau}
\affil[1]{ Universit\'e C\^ote d'Azur, CNRS, Inria, I3S, Sophia-Antipolis, France}
\affil[2]{LIRMM, Université de Montpellier, CNRS, France}
\date{}
\maketitle
\vspace{-10mm}
\begin{center}
{\small 
\texttt{$\{$nicolas.nisse, lucas.picasarri-arrieta$\}$@inria.fr}\\ 
\texttt{ignasi.sau@lirmm.fr}\\
}
\end{center}

\maketitle

\begin{abstract}

Given a digraph $D=(V,A)$ on $n$ vertices and a vertex $v\in V$, the \emph{cycle-degree} of $v$ is the minimum size of a set $S \subseteq V(D) \setminus \{v\}$ intersecting every directed cycle of $D$ containing $v$. From this definition of cycle-degree, we define the $c$-degeneracy (or cycle-degeneracy) of $D$, which we denote by $\delta^*_c(D)$. It appears to be a nice generalisation of the undirected degeneracy. For instance, the dichromatic number $\dic(D)$ of $D$ is bounded above by $\delta^*_c(D)+1$, where $\dic(D)$ is the minimum integer $k$ such that $D$ admits a $k$-dicolouring, i.e., a partition of its vertices into $k$ acyclic subdigraphs.

In this work, using this new definition of cycle-degeneracy, we extend several evidences for Cereceda's conjecture~\cite{cerecedaTHESIS} to digraphs. The \emph{$k$-dicolouring} graph of $D$, denoted by $\mathcal{D}_k(D)$, is the undirected graph whose vertices are the $k$-dicolourings of $D$ and in which two $k$-dicolourings are adjacent if they differ on the colour of exactly one vertex.  This is a generalisation of the $k$-colouring graph of an undirected graph $G$, in which the vertices are the proper $k$-colourings of $G$.

We show that $\mathcal{D}_k(D)$ has diameter at most $O_{\delta^*_c(D)}(n^{\delta^*_c(D) + 1})$ (respectively $O(n^2)$ and $(\delta^*_c(D)+1)n$) when $k$ is at least $\delta^*_c(D)+2$ (respectively $\frac{3}{2}(\delta^*_c(D)+1)$ and $2(\delta^*_c(D)+1)$). This improves known results on digraph redicolouring (Bousquet et al.~\cite{bousquetArXiv23}). 
Next, we extend a result due to Feghali~\cite{feghaliJCT147} to digraphs, showing that $\mathcal{D}_{d+1}(D)$ has diameter at most $O_{d,\epsilon}(n(\log n)^{d-1})$ when $D$ has maximum average cycle-degree at most $d-\epsilon$.

We then show that two proofs of Bonamy and Bousquet~\cite{bonamyEJC69} for undirected graphs can be extended to digraphs. The first one uses the digrundy number of a digraph $\dic_g(D)$, which is the worst number of colours used in a greedy dicolouring. If $k\geq \dic_g(D)+1$, we show that $\mathcal{D}_k(D)$ has diameter at most $4\cdot \dic(D) \cdot n$. The second one uses the \emph{$\mathscr{D}$-width} of a digraph, denoted by $\dw(D)$, which is a generalisation of the treewidth to digraphs. If $k\geq \dw(D)+2$, we show that $\mathcal{D}_k(D)$ has diameter at most $2(n^2+n)$.
 
Finally, we give a general theorem which makes a connection between the recolourability of a digraph $D$ and the recolourability of its underlying graph $\UG(D)$.
Assume that $\mathcal{G}$ is a class of undirected graphs, closed under edge-deletion and with bounded chromatic number, and let $k\geq \chi(\mathcal{G})$ (i.e., $k\geq \chi(G)$ for every $G\in \mathcal{G}$) be such that, for every $n$-vertex graph $G\in \mathcal{G}$, the diameter of the $k$-colouring graph of $G$ is bounded by $f(n)$ for some function $f$. We show that, for every $n$-vertex digraph $D$ such that $\UG(D)\in \mathcal{G}$, the diameter of $\mathcal{D}_k(D)$ is bounded by $2f(n)$. For instance, this result directly extends a number of results on planar graph recolouring to planar digraph redicolouring.

\end{abstract}

\bibliographystyle{plain}

\section{Introduction to graph recolouring}
\label{section:introduction}

 Given an undirected graph $G=(V,E)$ and a positive integer $k$, a \textit{$k$-colouring} of $G$ is a function $\alpha : V \xrightarrow{} [k] = \{1,\dots,k\}$. It is \emph{proper} if, for every edge $xy\in E$, we have $\alpha(x) \neq \alpha(y)$. So, for every $i\in [k]$, $\alpha^{-1}(i)$ induces an independent set on $G$. The \textit{chromatic number} of $G$, denoted by $\chi(G)$, is the smallest $k$ such that $G$ admits a proper $k$-colouring.  
  
 For any $k\geq \chi(G)$, the \textit{$k$-colouring graph} of $G$, denoted by ${\cal C}_k(G)$, is the graph whose vertices are the proper $k$-colourings of $G$ and in which two $k$-colourings are adjacent if they differ on the colour of exactly one vertex.
 A path between two given colourings in ${\cal C}_k(G)$ corresponds to a {\it sequence of recolourings}, that is, a sequence of pairs composed of a vertex of $G$, which is going to receive a new colour, and a new colour for this vertex.
 If ${\cal C}_k(G)$ is connected, we say that $G$ is \textit{$k$-mixing}. 
 In the last fifteen years, since the papers of Bonsma, Cereceda, van den Heuvel and Johnson~\cite{bonsmaTCS410,cerecedaEJC30,cerecedaJGT67}, graph recolouring has been studied by many researchers in graph theory.
 We refer the reader to the PhD thesis of Bartier~\cite{bartierTHESIS} for a complete overview on graph recolouring and to the surveys of van Heuvel~\cite{heuvel13} and Nishimura~\cite{Nishimura18} for reconfiguration problems in general.

A famous open conjecture on graph recolouring is due to Cereceda~\cite{cerecedaTHESIS} and makes a connection between the degeneracy of a graph and the diameter of ${\cal C}_k(G)$.
The \textit{degeneracy} of a graph $G$, denoted by $\delta^*(G)$, is the largest minimum degree of any subgraph of $G$.
Bonsma and Cereceda~\cite{bonsmaTCS410}  and  Dyer et al.~\cite{dyerRSA29} independently proved the following.
\begin{theorem}[Bonsma and Cereceda~\cite{bonsmaTCS410} ; Dyer et al.~\cite{dyerRSA29}]
    \label{thm:bonsma-cereda-dyer}
    Let $k\in \mathbb{N}$ and $G$ be a graph.
    If $k\geq \delta^*(G) +2$, then $G$ is $k$-mixing. 
\end{theorem}
The \textit{diameter} of $G$, denoted by $\diam(G)$, is the length of a longest shortest path of $G$.
The original proof of Theorem~\ref{thm:bonsma-cereda-dyer} also implies that $\mathcal{C}_k$ has diameter at most $2^n$, where $n=|V(G)|$.
Cereceda's conjecture states that the diameter of $\mathcal{C}_k$ is actually quadratic in $n$.
\begin{conjecture}[Cereceda~\cite{cerecedaTHESIS}]
    \label{conj:cereceda}
    Let $k\in \mathbb{N}$ and $G$ be a graph.
    If $k\geq \delta^*(G)+2$, then $\diam({\cal C}_k(G)) = O(n^2)$.
\end{conjecture}
In the remainder of this section, we recall several results approaching this conjecture, which can be seen as evidences for the general conjecture.
The following are the best existing bounds approaching Conjecture~\ref{conj:cereceda} in the general case\footnote{Given two computable functions $f,g$ and a parameter $\Gamma$, $f(n) = O_\Gamma(g(n))$ means that there exists a computable function $h$ such that $f(n) = O(h(\Gamma)\cdot g(n))$. Also $\Gamma$ can be a list of parameters $\Gamma_1,\dots,\Gamma_r$ in which case $f(n) = O_\Gamma(g(n))$ means that $f(n) = O(h(\Gamma_1,\dots,\Gamma_r)\cdot g(n))$.}.

\begin{theorem}[Bousquet and Heinrich~\cite{bousquetJCTB155}]
    \label{thm:bousquet-heinrich}
    Let $k\in \mathbb{N}$ and $G$ be a graph. Then ${\cal C}_k(G)$ has diameter at most:
    \begin{enumerate}[label=(\roman*)]
        \item \label{thm:bousquet-heinrich-statement-1} $O(n^2)$ if $k\geq \frac{3}{2}(\delta^*(G) +1)$,
        \item $O_{\epsilon}(n^{\lceil \frac{1}{\epsilon} \rceil})$ if $k\geq (1+\epsilon)(\delta^*(G) +1)$ for some $\epsilon > 0$, and
        \item \label{thm:bousquet-heinrich-statement-3} $O_d(n^{d+1})$ if $k\geq \delta^*(G) +2 = d+2$.
    \end{enumerate}
\end{theorem}
Bousquet and Perarnau~\cite{bousquetEJC52} also proved the following.
\begin{theorem}[Bousquet and Perarnau~\cite{bousquetEJC52}]
\label{thm:bousquet-perarnau}
    Let $k\in \mathbb{N}$ and $G$ be a graph.
    If $k\geq 2\delta^*(G)+2$, then $\diam({\cal C}_k(G)) \leq (\delta^*(G)+1)n$.
\end{theorem}

\medskip

In order to obtain possibly simpler versions of Conjecture~\ref{conj:cereceda}, one can restrict to graphs with bounded maximum average degree. The \textit{maximum average degree} of a graph $G$ is $\Mad(G) = \max \{ \frac{2|E(H)|}{|V(H)|} \mid H \text{ subgraph of } G\}$. It is easy to see that every graph $G$ satisfies $\lfloor \Mad(G) \rfloor \geq \delta^*(G)$. Hence Conjecture~\ref{conj:cereceda} would imply that every graph $G$ with $\Mad(G) \leq d-\epsilon$, where $d\in \mathbb{N}$ and $\epsilon > 0$, satisfies $\diam (\mathcal{C}_k(G)) = O(n^2)$ for every $k\geq d+1$. Feghali showed the following analogue result.
\begin{theorem}[Feghali~\cite{feghaliJCT147}]
    \label{thm:feghali}
    Let $d,k$ be integers such that $k \geq d+1$. For every $\epsilon > 0$ and every graph $G$ with $n$ vertices and maximum average degre at most $d-\epsilon$, $\diam(\mathcal{C}_k(G))=O_{d,\epsilon}(n(\log n)^{d-1})$.
\end{theorem}

Another approach toward Conjecture~\ref{conj:cereceda} consists of considering the maximum degree of a graph instead of its degeneracy. Cereceda has shown (see~\cite[Proposition~5.23]{cerecedaTHESIS}) that, for every graph $G$ on $n$ vertices and integer $k\geq \Delta(G)+2$, $\mathcal{C}_k(G)$ has diameter at most $(\Delta(G)+1)n$, where $\Delta(G)$ denotes the maximum degree of $G$. In order to get a more precise bound, Bonamy and Bousquet considered the grundy number. Let $G$ be a graph and $\mathcal{O}=(x_1,\dots,x_n)$ be an ordering  of $V(G)$. The \emph{greedy colouring} $\alpha_g(\mathcal{O},G)$ is the proper colouring in which every vertex $x_i$ receives the smallest colour that does not appear in $N(x_i) \cap \{x_1,\dots,x_{i-1} \}$. The \textit{grundy number} of $G$, denoted by $\chi_g(G)$, is the maximum, over all orderings $\mathcal{O}$, of the number of colours used in $\alpha_g(\mathcal{O},G)$. 

\begin{theorem}[Bonamy and Bousquet~\cite{bonamyEJC69}]
    \label{thm:grundy}
    For any graph $G$ on $n$ vertices, if $k \geq \chi_g(G) + 1$, then $G$ is $k$-mixing and $\diam(\mathcal{C}_k(G))\leq 4\cdot \chi(G) \cdot n$.
\end{theorem}

A last result approaching Conjecture~\ref{conj:cereceda} is due to Bonamy and Bousquet, and makes a connection between the treewidth of a graph and its recolourability. A \textit{tree-decomposition} of a graph $G=(V,E)$ is a pair $(T,\mathcal{X})$ where $T=(I,F)$ is a tree, and $\mathcal{X}=(B_i)_{i\in I}$ is a family of subsets of $V(G)$, called \textit{bags} and indexed by the vertices of $T$, such that:
\begin{enumerate}
    \item each vertex $v\in V$ appears in at least one bag, \textit{i.e.} $\bigcup_{i\in I} B_i= V$,
    \item for each edge $e = xy \in E$, there is an $i\in I$ such that $x,y \in B_i$, and 
    \item for each $v\in V$, the set of nodes indexed by $\{ i \mid i\in I, v\in B_i\}$ forms a subtree of $T$.
\end{enumerate}
The \textit{width} of a tree decomposition is defined as $\max_{i\in I} \{|B_i| -1\}$. The \textit{treewidth} of $G$, denoted by $\tw (G)$, is the minimum width of a tree-decomposition of $G$. It is easy to see that, for any graph $G$, $\tw (G) \geq \delta^*(G)$. Hence the following result is a weaker version of Conjecture~\ref{conj:cereceda}.
\begin{theorem}[Bonamy and Bousquet~\cite{bonamyEJC69}]
\label{thm:bonamy_bousquet_treewidth}
    Let $k\in \mathbb{N}$ and $G$ be an $n$-vertex graph.
    If $k\geq \tw(G) +2$, then $\diam({\cal C}_k(G)) \leq 2(n^2+n)$.
\end{theorem}

Bonamy et al.~\cite{bonamyJCO27} have shown that the bound of Theorem~\ref{thm:bonamy_bousquet_treewidth} is asymptotically sharp (up to a constant factor).

\paragraph{Organization of the paper.}  In this work, we extend every result mentioned above to digraphs and dicolourings. A $k$-dicolouring of a digraph $D$ is a partition of its vertices into at most $k$ acyclic subdigraphs. Every graph can be seen as a digraph in which every edge is actually two arcs, between the same vertices, in opposite directions. By extending a result, we mean showing that the result is actually true for every digraph and not only for symmetric digraphs. Note that a $k$-dicolouring of a symmetric digraph is indeed a proper $k$-colouring of its underlying graph.
More precisely, in Section~\ref{section:intro_digraphs}, we present the notations and the definitions used throughout the paper. In particular, we introduce the notion of cycle-degeneracy of a digraph which appears to be useful to play in directed graphs the role that the degeneracy plays in graphs. Then, we formally state our results, which are summarized in Table~\ref{table:bounded_diameter_digraphs}.

\renewcommand{\arraystretch}{1.5}
\begin{table}[H]
\begin{center}
\begin{tabular}{ |c || c | c | c | c | c |c|} \hline
    $k \geq $ & $ d+2$ & $ \frac{3}{2}(d+1)$ & $ 2(d+1)$  & $ \lceil \Mad_c \rceil + 1$ & $\dic_g+1$ & $ \dw+2$\\ \hline
    
    $\diam(\mathcal{D}_k(D)) $ & $O_d(n^{d+1})$ & $O(n^2)$ & $\leq (d+1)n$  & $O_{\Mad_c,\epsilon}(n(\log n)^{\lfloor \Mad_c \rfloor})$ & $\leq 4 \cdot \dic \cdot n$ & $O(n^2)$\\ \hline
    Theorem & \ref{thm:nd1} & \ref{thm:results_cycle_degeneracy}(ii) & \ref{thm:results_cycle_degeneracy}(iii) & \ref{thm:mad} & \ref{thm:digrundy}  & \ref{thm:dwidth}\\
    \hline
\end{tabular}

\caption{Bounds on the diameter of the $k$-dicolouring graph $\mathcal{D}_k(D)$ where $D$ is a digraph on $n$ vertices, with cycle-degeneracy $d$, maximum average cycle-degree $\Mad_c$, $\epsilon = \lceil \Mad_c \rceil - \Mad_c$, digrundy number $\dic_g$, dichromatic number $\dic$ and $\mathscr{D}$-width $\dw$. The notions above respectively extend the ones of $k$-colouring graph, degeneracy, maximum average degree, grundy number, chromatic number and treewidth to digraphs.}
\label{table:bounded_diameter_digraphs}
\end{center}
\end{table}

Section~\ref{section:cycle_degeneracy} is devoted to the proof of Theorem~\ref{thm:results_cycle_degeneracy}, which extends Theorem~\ref{thm:bonsma-cereda-dyer}, item (i) of Theorem~\ref{thm:bousquet-heinrich} and Theorem~\ref{thm:bousquet-perarnau} in the case of digraphs when the degeneracy is replaced by the cycle-degeneracy. In Section~\ref{section:nd1_mad}, we prove Theorems~\ref{thm:nd1} and~\ref{thm:mad}, which extend to digraphs item (iii) of Theorem~\ref{thm:bousquet-heinrich} and Theorem~\ref{thm:feghali} respectively. In Section~\ref{section:digrundy_number} we prove Theorem~\ref{thm:digrundy}, which generalises Theorem~\ref{thm:grundy} to digraphs.
In Section~\ref{section:general_bound_from_UG}, we prove Theorem~\ref{thm:general_bound_from_UG} which establishes a relationship between the diameter of the $k$-dicolouring graph of a digraph $D$ and the diameter of the $k$-colouring graph of its underlying graph $\UG(D)$, when $D$ belongs to a number of digraph classes (e.g., including planar digraphs).
Finally, in Section~\ref{section:directed_treewidth} we prove Theorem~\ref{thm:dwidth} which extends Theorem~\ref{thm:bonamy_bousquet_treewidth} to digraphs replacing the treewidth by the $\mathscr{D}$-width.
We conclude in Section~\ref{section:open_problems} by discussing the consequences of our results, especially on planar digraph redicolouring, and detail a few related open questions.

\section{Digraph redicolouring and main results}
\label{section:intro_digraphs}
    We refer the reader to~\cite{bang2009} for notation and terminology
not explicitly defined in this paper. Let $D=(V,A)$ be a digraph. A \textit{digon} is a pair of arcs in opposite directions between the same vertices. A \textit{simple arc} is an arc which is not in a digon. 
An \textit{oriented graph} is a digraph with no digons. The \textit{bidirected graph} associated with a graph $G$, denoted by $\bid{G}$,  is the digraph obtained from $G$ by replacing every edge by a digon. The \textit{underlying graph} of $D$, denoted by $\UG(D)$, is the undirected graph $G$ with vertex set $V(D)$ in which $uv$ is an edge if and only if $uv$ or $vu$ is an arc of $D$. 
Let $v$ be a vertex of a digraph $D$. The {\it out-degree} (resp. {\it in-degree}) of $v$, denoted by $d^+(v)$ (resp. $d^-(v)$), is the number of arcs leaving (resp. entering) $v$. We define the \textit{maximum degree} of $v$ as $d_{\max}(v) = \max\{d^+(v), d^-(v)\}$, and the \textit{minimum degree} of $v$ as $d_{\min}(v) = \min\{d^+(v), d^-(v)\}$. 
We also introduce the notion of cycle-degree. The \textit{cycle-degree} of $v$, denoted by $d_c(v)$, is the minimum size of a set $S \subseteq (V\setminus \{v\})$ such that $S$ intersects every directed cycle of $D$ containing $v$. 

\begin{sloppypar}
For each parameter $\Gamma \in \{\max, \min, c\}$ we define the corresponding \textit{maximum degree} of $D$ as $\Delta_\Gamma(D)=\max_{v\in V(D)} (d_{\Gamma}(v))$. Analogously, we define the corresponding \textit{minimum degree} of $D$ as $\delta_\Gamma(D) = \min_{v\in V(D)} (d_{\Gamma}(v))$. Finally, we define the \textit{$\Gamma$-degeneracy}  of $D$, denoted by $\delta^*_\Gamma(D)$, as $\max_{H\subseteq D}(\delta_\Gamma(H))$, where $H\subseteq D$ denotes that $H$ is a (not necessarily induced) subdigraph of $D$. An equivalent characterisation of the degeneracy is the following. For each parameter $\Gamma$, $\delta^*_\Gamma(D) \leq k$ if and only if there exists an ordering $v_1,\dots,v_n$ of $V(D)$ such that, for each $i\in [n]$, then $d_{\Gamma}(v_i) \leq k$ in $D[\{v_i,\dots,v_n\}]$.  
\end{sloppypar}

Observe that we have $d_c(v) \leq d_{\min}(v) \leq d_{\max}(v)$. The first inequality holds because each of the sets $N^+(v)$ and $N^-(v)$ intersects every directed cycle of $D$ containing $v$, and the second inequality holds by definition. This implies $\Delta_c(D) \leq \Delta_{\min}(D) \leq \Delta_{\max}(D)$ and $ \delta^*_{c}(D) \leq \delta^*_{\min}(D) \leq \delta^*_{\max}(D)$. Moreover, if $N^+(v) = N^-(v)$ (meaning that $v$ is not incident to a simple arc) then we have $d_c(v) = d_{\min}(v) = d_{\max}(v)$. Hence, if $D$ is a bidirected graph $\bid{G}$, for each parameter $\Gamma \in \{c,\min,\max\}$, we have $\Delta(G) = \Delta_\Gamma(\bid{G})$ and $\delta^*(G) = \delta^*_{\Gamma}(\bid{G})$.
As we will show later, the notion of $c$-degeneracy appears to be a natural generalisation of the undirected degeneracy when dealing with the directed treewidth. To the best of the authors' knowledge, this extension of the classical degeneracy has not been considered before.

In 1982, Neumann-Lara~\cite{neumannlaraJCT33} introduced the notions of dicolouring and dichromatic number, which generalize the ones of proper colouring and chromatic number. 
For a positive integer $k$, a \textit{$k$-colouring} of $D=(V,A)$ is a function $\alpha : V \xrightarrow{} [k]$. It is a \textit{$k$-dicolouring} if $\alpha^{-1}(i)$ induces an acyclic subdigraph in $D$ for each $i \in [k]$. The \textit{dichromatic number} of $D$, denoted by $\dic(D)$, is the smallest $k$ such that $D$ admits a $k$-dicolouring. 
There is a one-to-one correspondence between the proper $k$-colourings of a graph $G$ and the $k$-dicolourings of its associated bidirected graph $\bid{G}$, and in particular $\chi(G) = \dic(\bid{G})$. Hence every result on graph proper colourings can be seen as a result on dicolourings of bidirected graphs, and it is natural to study whether the result can be extended to all digraphs. 

Let us show that $\dic(D) \leq \delta^*_c(D) + 1$. By definition of $c$-degeneracy, we can find an ordering $v_1,\dots,v_n$ of $V(D)$ such that, for every $i\in [n]$, there exists $S_i \subseteq \{v_{i+1},\dots,v_n \}$  of size at most $\delta^*_c(D)$, and such that $S_i \cup \{v_{1},\dots,v_{i-1}\}$ intersects every directed cycle of $D$ containing $v_i$. Hence, considering the vertices from $v_n$ to $v_1$, one can greedily find a $(\delta_c^*(D)+1)$-dicolouring of $D$ by colouring each $v_i$ with a colour that has not been chosen in $S_i$. Indeed, suppose for contradiction that, for some $i \in [\delta^*_c(D) + 1]$, the subdigraph of $D$ induced by the set of vertices assigned colour $i$ contains a directed cycle $C$. Let $v_j$ be the leftmost vertex in $C$ according to the considered ordering. Then, by definition of  $c$-degeneracy we have that $(V(C) \setminus \{v_j\}) \cap  S_j\neq \emptyset$, but this contradicts the fact that, by construction of the colouring, the colour of $v_j$ is different from all the colours of the vertices in $S_j$.

The notion of digraph redicolouring was first introduced by Bousquet et al. in~\cite{bousquetArXiv23}.
For any $k\geq \dic(D)$, the \textit{$k$-dicolouring graph} of $D$, denoted by ${\cal D}_k(D)$, is the undirected graph whose vertices are the $k$-dicolourings of $D$ and in which two $k$-dicolourings are adjacent if they differ by the colour of exactly one vertex.  Observe that ${\cal C}_k(G) = {\cal D}_k(\bid{G})$ for any graph $G$.  A \textit{redicolouring sequence} between two dicolourings is a path between these dicolourings in ${\cal D}_k(D)$.  
The digraph $D$ is \textit{$k$-mixing} if ${\cal D}_k(D)$ is connected. 
In~\cite{bousquetArXiv23}, the authors mainly study the $k$-dicolouring graph of digraphs with bounded min-degeneracy or bounded maximum average degree, and they show that finding a redicolouring sequence between two given $k$-dicolouring of a digraph is PSPACE-complete. Dealing with the min-degeneracy of digraphs, they extended Theorems~\ref{thm:bonsma-cereda-dyer},~\ref{thm:bousquet-heinrich}\ref{thm:bousquet-heinrich-statement-1} and~\ref{thm:bousquet-perarnau}:
\begin{theorem}[Bousquet et al.~\cite{bousquetArXiv23}]
\label{thm:results_min_degeneracy}
    Let $k\in \mathbb{N}$ and $D$ be a digraph on $n$ vertices. Then:
    \begin{enumerate}[label=(\roman*)]
        \item if $k\geq \delta^*_{\min}(D)+2$ then $D$ is $k$-mixing,
        \item if $k\geq \frac{3}{2}(\delta^*_{\min}(D)+1)$ then $\diam( \mathcal{D}_k(D)) = O(n^2)$, and
        \item if $k\geq 2(\delta^*_{\min}(D) +1)$ then $\diam( \mathcal{D}_k(D)) \leq (\delta^*_{\min}(D)+1)n$.
    \end{enumerate}
\end{theorem}

They also proposed a stronger version of Conjecture~\ref{conj:cereceda} for digraphs:
\begin{conjecture}[Bousquet et al.~\cite{bousquetArXiv23}]
    \label{conj:directed_cereceda} 
    Let $k\in \mathbb{N}$ and $D$ be a digraph. If $k\geq \delta^*_{\min}(D) + 2$, then $\diam(\mathcal{D}_k(D)) = O(n^2)$.
\end{conjecture}

Digraph redicolouring was also investigated in~\cite{picasarriArxiv23}, where the author deals with the maximum degrees of a digraph instead of its degeneracies. In particular, he shows the following extension of a result of Feghali et al.~\cite{feghaliJGT83} to digraphs. We say that a digraph is {\it connected} if its underlying graph is connected.
\begin{theorem}[Picasarri-Arrieta~\cite{picasarriArxiv23}]
    \label{thm:k_delta_1}
    Let $D=(V,A)$ be a connected digraph with $\Delta_{\max}(D) = \Delta \geq 3$, $k\geq \Delta+1$, and $\alpha$, $\beta$ be two $k$-dicolourings of $D$. Then one of the following holds:
    \begin{itemize}
        \item $\alpha$ or $\beta$ is an isolated vertex in $\mathcal{D}_k(D)$, or
        \item there is a redicolouring sequence of length at most $O(\Delta^2|V|^2)$ between $\alpha$ and $\beta$.
    \end{itemize}
\end{theorem}

In Section~\ref{section:cycle_degeneracy}, we show that each of the proofs of the statements of Theorem~\ref{thm:results_min_degeneracy} can be adapted to show the following stronger result.
\begin{restatable}{theorem}{resultscycledegeneracy}
    \label{thm:results_cycle_degeneracy}
    Let $k\in \mathbb{N}$ and $D$ be a digraph on $n$ vertices. Then:
    \begin{enumerate}[label=(\roman*)]
        \item \label{thm:results_cycle_degeneracy_statement_1} if $k\geq \delta^*_{c}(D)+2$ then $D$ is $k$-mixing, 
        \item \label{thm:results_cycle_degeneracy_statement_2} if $k\geq \frac{3}{2}(\delta^*_{c}(D)+1)$ then $\diam( \mathcal{D}_k(D)) = O(n^2)$, and
        \item \label{thm:results_cycle_degeneracy_statement_3} if $k\geq 2(\delta^*_{c}(D) +1)$ then $\diam( \mathcal{D}_k(D)) \leq (\delta^*_{c}(D)+1)n$.
    \end{enumerate}
\end{restatable}

The \textit{maximum average degree} of $D$ is $\Mad(D) = \max \{ \frac{2|A(H)|}{|V(H)|} \mid H \text{ subdigraph of } D\}$. The average cycle-degree of a digraph $D=(V,A)$ is $\Ad_c(D) = \frac{1}{|V|}\sum_{v\in V}d_c(v)$. The \textit{maximum average cycle-degree} of $D$ is $\Mad_c(D) = \max \{ \Ad_c(H) \mid H \text{ subdigraph of } D\}$. For every undirected graph $G$, we have $\Mad_c(\bid{G}) = \Mad(G)$.

In Section~\ref{section:nd1_mad}, we prove, based on the proofs of~\cite{feghaliJCT147}, the following two extensions of Theorem~\ref{thm:bousquet-heinrich}\ref{thm:bousquet-heinrich-statement-3} and Theorem~\ref{thm:feghali}, respectively. 
\begin{restatable}{theorem}{ndone}
    \label{thm:nd1}
    Let $D$ be a digraph and $k\geq \delta^*_c(D) +2 = d+2$. Then $\diam(\mathcal{D}_k(D)) = O_d(n^{d+1})$.
\end{restatable}

\begin{restatable}{theorem}{mad}
    \label{thm:mad} 
    Let $d \geq 1$ and $k\geq d+1$ be two integers, and let $\epsilon > 0$. If $D$ is a digraph satisfying $\Mad_{c}(D)\leq d - \epsilon$, then $\diam(\mathcal{D}_{k}(D)) = O_{d,\epsilon}(n(\log n)^{d-1})$.
\end{restatable}

 Observe that every digraph $D$ satisfies $\Mad_c(D) \leq \frac{1}{2}\Mad(D)$. This holds because, for every vertex $v\in V(D)$, $d_c(v) \leq \frac{1}{2}(d^+(v) + d^-(v))$. Hence, for every subdigraph $H$ of $D$, we have:
 \[ 2|A(H)| = \sum_{v\in V(H)}(d^+(v) + d^-(v)) \geq 2\sum_{v\in V(H)}d_c(v) = 2\cdot \Ad_c(H) \cdot |V(H)|.\]
Hence, the following is a direct consequence of Theorem~\ref{thm:mad}.
\begin{corollary}
    \label{cor:mad} 
    Let $d \geq 1$,$k\geq \lfloor \frac{d+3}{2} \rfloor$ be two integers, and let $\epsilon > 0$. If $D$ is a digraph satisfying $\Mad(D)\leq d - \epsilon$, then $\mathcal{D}_{k}(D)$ has diameter at most:
    \begin{enumerate}[label=(\roman*)]
        \item $O_{d}\left(n(\log n)^{\frac{d-1}{2}}\right)$ if $d$ is odd, and
        \item $O_{d,\epsilon}\left(n(\log n)^{\frac{d-2}{2}}\right)$ otherwise.
    \end{enumerate}
\end{corollary}

The \textit{digrundy number} of a digraph $D=(V,A)$, introduced in~\cite{araujoDAM317}, is the natural analogue of the grundy number for digraphs. If $\phi$ is a dicolouring of $D$, then $\phi$ is a \textit{greedy dicolouring} if there is an ordering $v_1,\dots,v_n$ of $V$ such that, for each vertex $v_i$ and each colour $c$ smaller than $\phi(v_i)$, the set of vertices $(\{v_1,\dots,v_{i-1} \} \cap \phi^{-1}(c)) \cup \{v_i\}$ contains a directed cycle. The digrundy number of $D$, denoted by $\dic_g(D)$, corresponds to the maximum number of colours used in a greedy dicolouring of $D$. The following result generalises Theorem~\ref{thm:grundy}.

\begin{restatable}{theorem}{digrundy}
    \label{thm:digrundy}
    For any digraph $D$, if $k\geq \dic_g(D)+1$, then $\diam( \mathcal{D}_k(D)) \leq 4 \cdot \dic(D) \cdot n$.
\end{restatable}
Analogously to the undirected case, we always have $\dic(D) \leq \dic_g(D) \leq \Delta_c(D) + 1$. Thus, the following is a direct consequence of Theorem~\ref{thm:digrundy}. Note also that this is an improvement of the case $k\geq \Delta_{\max}(D) + 2$ of Theorem~\ref{thm:k_delta_1}.
\begin{corollary}
    \label{thm:k_delta_2}
    For any digraph $D$, if $k\geq \Delta_c(D) +2$, then $\diam(\mathcal{D}_k(D)) \leq 4\cdot \dic(D) \cdot n$.
\end{corollary}

In Section~\ref{section:general_bound_from_UG} we show the following general result which makes a connection between the recolourability of a digraph and the recolourability of its underlying graph. 
\begin{restatable}{theorem}{generalbound}
    \label{thm:general_bound_from_UG}
    Let $\mathcal{G}$ be a family of undirected graphs, closed under edge-deletion and with bounded chromatic number, and let $k \geq \chi(\mathcal{G})$ (\textit{i.e.} $k \geq \chi(G)$ for every $G\in \mathcal{G}$) be such that, for every $n$-vertex graph $G\in \mathcal{G}$, the diameter of $\mathcal{C}_k(G)$ is bounded by $f(n)$ for some function $f$. Then for any $n$-vertex digraph $D$ such that $\UG(D) \in \mathcal{G}$, the diameter of $\mathcal{D}_k(D)$ is bounded by $2f(n)$. 
\end{restatable}
In particular, since removing edges does not increase the treewidth of a graph, the following is a consequence of Theorems~\ref{thm:bonamy_bousquet_treewidth} and~\ref{thm:general_bound_from_UG} (by taking $\mathcal{G}= \{G \mid \tw(G) \leq \ell \}$ for some constant $\ell$).
\begin{corollary}
    Let $k\in \mathbb{N}$ and $D$ be a digraph.
    If $k\geq \tw(\UG(D)) +2$, then $\diam({\cal D}_k(D)) = O(n^2)$.
\end{corollary}
However, the treewidth of the underlying graph of $D$ is not a satisfying extension of  treewidth to digraphs, since it does not take under consideration the orientations in $D$. There exist, at least, four well known generalisations of  treewidth to digraphs: the directed treewidth (introduced in~\cite{johnsonJCT82}, see also~\cite{reedENDM3}), the $\mathscr{D}$-width (introduced in~\cite{safariLNCS3618}, see also~\cite{safariThesis}), the DAG-width (introduced in~\cite{berwangerJCTB102}) and the Kelly-width (introduced in~\cite{hunterTCS399}).

An \emph{out-arborescence} is a rooted tree in which every edge is oriented away from the root. A {\it directed tree-decomposition} $(T,{\cal W},{\cal X})$ of a digraph $D=(V,A)$ consists of an out-arborescence $T=(I,F)$ rooted in $r\in I$, a partition ${\cal W}=(W_t)_{t\in I}$ of $V$ into non-empty parts, and a family  ${\cal X}=(X_e)_{e\in F}$ of subsets of vertices of $D$ such that, for every $tt'\in F$ we have:
\begin{enumerate}
    \item $X_{tt'}\cap (\bigcup_{t'' \in T_{t'}} W_{t''}) = \emptyset$ (where $T_{t'}$ denotes the subtree of $T$ rooted in $t'$), and
    \item for every directed walk $P$ with both ends in $\bigcup_{t'' \in T_{t'}} W_{t''}$ and some internal vertex not in $\bigcup_{t'' \in T_{t'}} W_{t''}$, it holds that $V(P) \cap X_{tt'} \neq \emptyset$.
\end{enumerate}

The {\it width} of $(T,{\cal W},{\cal X})$ equals $\max_{t \in I} |H_t|-1$, where $H_t= W_t \cup \bigcup_{e \in F, t\in e}X_{e}$, and the {\it directed treewidth} of $D$, denoted by $\dtw(D)$, is the minimum width of its directed tree-decompositions. 
Recall that the treewidth of an undirected graph is always at least its degeneracy. However it is well known that there exist digraphs with arbitrary large min-degeneracy and directed tree-width exactly one. We include a proof for completeness.
\begin{proposition}[Folklore]
    \label{prop:dtw_1_min_d}
    For every integer $d$, there exists a digraph $D=(V,A)$ such that every vertex $v\in V$ satisfies $d^+(v) \geq d$, $d^-(v) \geq d$, and $\dtw(D) = 1$.
\end{proposition}
\begin{proof}
    Let $T$ be a tree rooted in $r\in V(T)$ with depth at least $d$ (that is, all leaves are at distance at least $d$ from the root), such that every non-leaf vertex has at least $d$ children. We orient each edge $uv$ of $T$ from the parent to its child. Then we add every arc $uv$ such that $v$ is an ancestor of $u$. In the obtained digraph $D$, every vertex has out-degree at least $d$. 

    Then we add a disjoint copy $\cev{D}$ of $D$ in which we reverse every arc, so in $\cev{D}$ every vertex has in-degree at least $d$. We finally add every arc from $\cev{D}$ to $D$.

    In the resulting digraph, every vertex has out-degree and in-degree at least $d$. Moreover, the directed treewidth of a digraph is equal to the maximum directed treewidth of its strongly connected components.
    For each edge $uv$ of $T$, such that $u$ is the parent of $v$, we label $uv$ with $u$. Then $T$, together with this labelling, is a directed tree-decomposition of both $D$ and $\cev{D}$ and has width 1. This follows from the fact that every directed cycle of $D$ containing a vertex must also contain its father in $T$.
\end{proof}

\begin{sloppypar}
The following proposition shows that, dealing with directed treewidth,  $c$-degeneracy, compared to min-degeneracy, appears to be a better generalisation of the undirected one.
\end{sloppypar}
\begin{proposition}
    \label{prop:dtw_c_deg}
    For every digraph $D$, it holds that $\dtw(D) \geq \delta^*_c(D)$.
\end{proposition}
\begin{proof}
Consider an optimal directed tree-decomposition $(T,\mathcal{W},\mathcal{X})$ of $D$. If $t$ is a leaf of $T$ and $v$ is a vertex in $W_t$, then $H_t\setminus \{v\}$ intersects every directed cycle of $D$ containing $v$. Thus, $d_c(v) \leq \dtw(D)$. Moreover, since $t$ is a leaf of $T$, it is easy to verify that removing $v$ from $D$ does not increase its directed treewidth, and we can repeat the same argument in $D \setminus \{v\}$.
\end{proof}

Proposition~\ref{prop:dtw_c_deg} implies that $\dic(D) \leq \dtw(D) + 1$ for every digraph $D$ and, together with Theorems~\ref{thm:results_cycle_degeneracy} and~\ref{thm:nd1}, implies the following.
\begin{corollary}
    \label{cor:dtw_recolourability}
    Let $k\in \mathbb{N}$ and $D$ be a digraph on $n$ vertices.
    \begin{enumerate}[label=(\roman*)]
        \item if $k\geq \dtw(D)+2$ then $D$ is $k$-mixing and $\diam(\mathcal{D}_k(D)) = O_{\dtw}(n^{\dtw+1})$,
        \item if $k\geq \frac{3}{2}(\dtw(D)+1)$ then $\diam( \mathcal{D}_k(D)) = O(n^2)$, and
        \item if $k\geq 2(\dtw(D) +1)$ then $\diam( \mathcal{D}_k(D)) \leq (\dtw(D)+1)n$.
    \end{enumerate}
\end{corollary}
The interested reader may have a look at~\cite{berwangerJCTB102} (respectively~\cite{hunterTCS399}) and see that DAG-width (respectively Kelly-width) is bounded below by min-degeneracy (respectively by min-degeneracy plus one). Thus, Corollary~\ref{cor:dtw_recolourability} also holds for DAG-width (respectively Kelly-width minus one).

\medskip

Finally, in Section~\ref{section:directed_treewidth}, we show that the proof of Theorem~\ref{thm:bonamy_bousquet_treewidth} extends to digraphs using $\mathscr{D}$-width. A {\it $\mathscr{D}$-decomposition} of a digraph $D=(V,A)$ is a pair $(T,{\cal X})$ such that $T=(I,F)$ is an undirected tree and $\mathcal{X}=(X_v)_{v \in I}$ is a family of subsets (called {\it bags}) of $V$ indexed by the nodes of $T$, which satisfies Property~$(*)$ stated below, for which we first need a definition. For a vertex subset $S \subseteq V$, let the {\it support} of $S$ in $(T,{\cal X})$, denoted by $T_S$, be the subgraph of $T$ with vertex-set $\{t \in I \mid X_t \cap S \neq \emptyset\}$ and edge-set $\{ tt' \in F \mid \exists u \in S \cap X_t \cap X_{t'}\}$. A $\mathscr{D}$-decomposition must ensure the following property: 

\medskip

$(*)$ for every subset $S \subseteq V$, if $D[S]$ is strongly connected then $T_S$ is a non-empty subtree of $T$. 

\medskip

 Similarly to undirected tree-decompositions, the {\it width} of $(T,{\cal X})$ is the maximum size of its bags minus one, and the {\it $\mathscr{D}$-width} of $D$, denoted by $\dw(D)$, is the minimum width of its $\mathscr{D}$-decompositions.
 Note that, for every $v \in V$, $D[\{v\}]$ is strongly connected, and therefore Property $(*)$ can be seen as a generalization of the basic properties of undirected tree-decompositions that state that the set of bags containing some vertex $v$ must induce a (connected) subtree and that every vertex must belong to at least one bag. Note also that, if $\{u,v\}$ is a digon of $D$, Property $(*)$ implies that $u$ and $v$ must belong to a common bag of $(T,\mathcal{X})$. Hence, every bidirected graph $G$ satisfies $\tw(G) = \dw(\bid{G})$, and the following is actually a generalisation of Theorem~\ref{thm:bonamy_bousquet_treewidth}. Our proof is strongly based on the proof of Theorem~\ref{thm:bonamy_bousquet_treewidth}.

\begin{restatable}{theorem}{dwidth}
    \label{thm:dwidth}
    If $D=(V,A)$ is an $n$-vertex digraph with $\dw(D)\leq k-1$, then $\diam(\mathcal{D}_{k+1}(D)) \leq 2(n^2+n)$.
\end{restatable}

Note that the bound of Theorem~\ref{thm:dwidth} is asymptotically sharp (up to a constant factor) since Theorem~\ref{thm:bonamy_bousquet_treewidth} is already known to be sharp.
Finally, observe that the digraph $D$ built in the proof of Proposition~\ref{prop:dtw_1_min_d} also satisfies $\dw(D) = 1$ but $\delta^*_{\min}(D) \geq d$. Again, the following easy proposition shows that, dealing with the $\mathscr{D}$-width, $c$-degeneracy, compared to min-degeneracy, appears to be a better generalisation of the undirected one.

\begin{proposition}
    \label{prop:dw_c_deg}
    For every digraph $D$, it holds that $\dw(D) \geq \delta^*_c(D)$.
\end{proposition}
\begin{proof}
    Consider an optimal directed $\mathscr{D}$-decomposition $(T,\mathcal{X}=(X_t)_{t\in V(T)})$ of $D$. Let $t$ be a leaf of $T$ and $v$ be a vertex in $X_t$ that belongs to no other bag $X_{t'}$ (this is possible unless $X_t \subseteq X_{t'}$, $tt'\in E(T)$, in which case we just remove the bag $t$ from the decomposition). We claim that $X_t\setminus \{v\}$ intersects every directed cycle of $D$ containing $v$ (which directly implies $d_c(v) \leq \dw(D)$). Assume not, and let $C$ be a directed cycle such that $X_t \cap V(C) = \{v\}$. Then, since $D[V(C)]$ is strongly connected, $T_{V(C)}$ must be connected. This is a contradiction since $t$ is an isolated vertex in  $T_{V(C)}$.
\end{proof}

Analogously to the directed treewidth, Proposition~\ref{prop:dw_c_deg}, together with Theorems~\ref{thm:results_cycle_degeneracy} and~\ref{thm:nd1}, implies the following (note that the two first items are also implied by Theorem~\ref{thm:dwidth}, but the third one is not).
\begin{corollary}
    \label{cor:dw_recolourability}
    Let $k\in \mathbb{N}$ and $D$ be a digraph on $n$ vertices.
    \begin{enumerate}[label=(\roman*)]
        \item if $k\geq \dw(D)+2$ then $D$ is $k$-mixing and $\diam(\mathcal{D}_k(D)) = O_{\dw}(n^{\dw+1})$,
        \item if $k\geq \frac{3}{2}(\dw(D)+1)$ then $\diam( \mathcal{D}_k(D)) = O(n^2)$, and
        \item if $k\geq 2(\dw(D) +1)$ then $\diam( \mathcal{D}_k(D)) \leq (\dw(D)+1)n$.
    \end{enumerate}
\end{corollary}

\section{Bounds on the diameter of \texorpdfstring{$\mathcal{D}_k(D)$}{Dk(D)} when \texorpdfstring{$k\geq \frac{3}{2}(\delta^*_c(D)+1)$}{k >= 3/2(dc*(D)+1)}}
\label{section:cycle_degeneracy}

This section is devoted to the proof of Theorem~\ref{thm:results_cycle_degeneracy}, and we will prove each item of its statement independently. Let us restate it.

\resultscycledegeneracy*

\begin{proof}[Proof of Theorem~\ref{thm:results_cycle_degeneracy}\ref{thm:results_cycle_degeneracy_statement_1}]
The proof is by induction on $n=|V(D)|$. The result clearly holds for $n\leq 1$. Let us assume that $n>1$ and that the result holds for $n-1$. 
Let $\alpha,\beta$ be any two $k$-dicolourings of $D$ and let $v \in V$ be a vertex satisfying $d_c(v) \leq \delta^*_c(D)$. Let $\alpha',\beta'$ be the two $k$-dicolourings induced, respectively, by $\alpha$ and $\beta$ on $D - \{v\}$. By induction, there exists a redicolouring sequence $\alpha'=\alpha'_1,\dots,\alpha'_q=\beta'$ where $\alpha'_i$ and $\alpha'_{i-1}$ differ by the colour of exactly one vertex $v_i \in V\setminus \{v\}$, for every $1<i\leq q$.

Now, we build the following redicolouring sequence from $\alpha$ to $\beta$. At step $i$, if $v_i$ can be recoloured as from $\alpha'_{i-1}$ to $\alpha'_{i}$, then recolour $v_i$ accordingly. Otherwise, this implies that there exists a directed cycle containing $v$ and $v_i$ whose all vertices (but $v_i$) have colour $\alpha'_i(v_i)$. By definition of  cycle-degree, there exists a transversal $X$ of the directed cycles containing $v$, with  $|X|\leq \delta^*_c(D)$ and $v \notin X$. Let $c \neq \alpha'_i(v_i)$ be a colour that does not appear in $X$ (it exists since $k\geq \delta^*_c(D)+2$). Colour $v$ with $c$ and then $v_i$ with $\alpha'_i(v_i)$. 
Finally (after step $q$), recolour $v$ with its final colour $\beta(v)$.
\end{proof}

A {\it $k$-list assignment} in an undirected graph $G$ is a function $L: V(G) \to 2^{[k]}$ on $V(G)$ such that  $L(v)\subseteq [k]$ for each vertex $v \in V(G)$. An {\it $L$-colouring} of $G$ is a proper colouring $\alpha$ of $G$ such that $\alpha(v) \in L(v)$ for every vertex $v$. Given such a list assignment $L$, ${\cal C}(G,L)$ is the graph whose vertices are the $L$-colourings of $G$ and in which two colourings are adjacent if they differ by the colour of exactly one vertex. An $L$-recolouring sequence is a path in ${\cal C}(G,L)$. Also, we say that $L$
is {\it $a$-feasible} if, for some ordering $v_1,\dots,v_n$ of $V$, $|L(v_i)| \geq |N(v) \cap \{v_{i+1},\dots,v_n\}| + 1 + a$ for every $i\in [n]$. 
Also we say that an $L$-colouring $c$ \emph{avoids a set of colours $S$} if for every vertex $v\in V(G)$, $c(v)$ does not belong to $S$. We need the following lemma from~\cite[Lemma~20]{bousquetArXiv23}, which was indirectly proved first in~\cite{bousquetJCTB155}.

\begin{lemma}[\!\!\cite{bousquetArXiv23, bousquetJCTB155}]
    \label{lemma:avoid-colours}
    Let $G=(V,E)$ be an undirected graph on $n$ vertices, $L$ be a $k$-list assignment of $G$ that is $\left \lceil \frac{k}{3} \right \rceil$-feasible and $\alpha$ an $L$-colouring of $G$ that avoids a set $S$ of $\left \lceil \frac{k}{3} \right \rceil$ colours. Then for any set of $\left \lceil \frac{k}{3} \right \rceil$ colours $S'$, there is an $L$-colouring $\beta$ of $G$ that avoids $S'$ and such that there is an $L$-recolouring sequence from $\alpha$ to $\beta$ of length at most $\frac{4k + 12}{3}n$.
\end{lemma}

\begin{proof}[Proof of Theorem~\ref{thm:results_cycle_degeneracy}\ref{thm:results_cycle_degeneracy_statement_2}]
    Let $D=(V,A)$ be a digraph on $n$ vertices and $k\geq \frac{3}{2}(\delta^*_{c}(D)+1)$. 
    Let $(v_1,\dots,v_n)$ be a cycle-degeneracy-ordering of $D$, that is, an ordering such that for each $i\in [n]$, there exists $X_i \subseteq \{v_{i+1},\dots,v_n\}$, $|X_i| \leq \delta^*_{c}(D)$, such that every directed cycle of $D$ containing $v_i$ must intersect $\{v_1,\dots,v_{i-1}\} \cup X_i$. 

    Let $G=(V,E)$ be the undirected graph where $E = \{ v_iv_j \mid v_j \in X_i, i\in [n] \}$. We first prove that each proper colouring of $G$ is a dicolouring of $D$. Assume that this is not the case, and there exists a proper colouring $\alpha$ of $G$ such that $D$, coloured with $\alpha$, contains a monochromatic directed cycle $C$. Let $v_i$ be the least vertex of $C$ in the ordering $(v_1,\dots,v_n)$. Then $C$ must contain a vertex $v_j$ in $X_i$. This is a contradiction, since $\alpha$ is a proper colouring of $G$ and $v_iv_j \in E$.
    
    By construction, $G$ has degeneracy at most $\delta^*_{c}(D)$. Using Theorem~\ref{thm:bousquet-heinrich}(i), we get that ${\cal C}_k(G)$ has diameter at most $C_0n^2$ for some constant $C_0$.
    
    \medskip
    
     Let $\alpha$ be any $k$-dicolouring of $D$. We will now show that there exists a dicolouring $\alpha'$ of $D$ that is also a proper colouring of $G$, and such that there exists a redicolouring sequence between $\alpha$ and $\alpha'$ of length at most $C_1n^2$ for some constant $C_1$.
     Set $\delta^*=\delta^*_{c}(D) \geq \delta^*(G)$, $Y_i=\{v_{i+1},\dots,v_n\}$, and $H_i = G-Y_i$ for all $i\in [n]$.

    Let $L_i$ be the $k$-list assignment of $H_i$ defined by
    $$L_i(v_j) = [k] \setminus \{ \alpha(v) \mid v\in X_j \cap Y_i \}~\mbox{for~all~} j\in [i].$$
      Since $k$, the total number of colours, is at least $\frac{3}{2}(\delta^* + 1)$, for every $j\in [i]$ we have:
    \begin{align*}
        |L_i(v_j)| &\geq k - |X_j \cap Y_i|\\
                    &\geq \frac{k}{3} + \frac{2}{3}\frac{3}{2}(\delta^* +1) - |X_j\cap Y_i|\\
                    &\geq |X_j\cap \{v_{j+1},\dots,v_i\}| + 1 + \frac{k}{3}.
    \end{align*}
    Hence, since $|L_i(v_j)|$ is an integer, $L_i$ is a $\left \lceil \frac{k}{3} \right \rceil$-feasible $k$-list assignment of $H_i$.
    
    \begin{remark}\label{remark:extendstoD}
    Let $\gamma$ be a dicolouring of $D$ such that for some $i$, $\gamma$ agrees with $\alpha$ on $\{v_{i+1},\dots,v_n\}$ and $\gamma_{|H_i}$ (the restriction of $\gamma$ to $H_i$) is an $L_i$-colouring of $H_i$. Then any $L_i$-recolouring sequence starting from $\gamma_{|H_i}$ on $H_i$ is a redicolouring sequence in $D$. Indeed, assume this is not the case and at one step, we get to an $L_i$-colouring $\zeta$ of $H_i$ but $\zeta_D$ contains a monochromatic cycle $C$, where $\zeta_D(v) = \zeta(v)$ when $v$ belongs to $H_i$ and $\zeta_D(v) = \gamma(v)$ otherwise. Let $v_j$ be the vertex of $C$ such that $j$ is minimum in the cycle-degeneracy-ordering of $D$. Then $C$ must intersect $X_j$ in some vertex $v_q$. Thus either $q\leq i$ and then $v_{q}v_j$ is a monochromatic edge in $H_i$ or $q\geq i+1$ but then $\zeta(v_q) = \gamma(v_q) = \alpha(v_q)$ does not belong to $L_i(v_j)$. In both cases, we get a contradiction.
    \end{remark}
    
    \begin{claim}\label{claim+k/3}
    Let $\gamma_i$ be a $k$-dicolouring of $D$, agreeing with $\alpha$ on $Y_{i}$, which induces an $L_i$-colouring of $H_i$ avoiding at least $\left \lceil \frac{k}{3} \right \rceil$ colours in $H_i$. 
    Then there is a redicolouring sequence of length at most $\frac{8k + 24}{3}n + \left \lceil \frac{k}{3} \right \rceil$ from $\gamma_i$ to a dicolouring $\gamma_{i+\left \lceil \frac{k}{3} \right \rceil}$ which induces an $L_{i+\left \lceil \frac{k}{3} \right \rceil}$-colouring of $H_{i+\left \lceil \frac{k}{3} \right \rceil}$ avoiding at least $\left \lceil \frac{k}{3} \right \rceil$ colours in $H_{i+\left \lceil \frac{k}{3} \right \rceil}$. Moreover, $\gamma_{i+\left \lceil \frac{k}{3} \right \rceil}$ agrees with $\alpha$ on $Y_{i+\left \lceil \frac{k}{3} \right \rceil}$.
    \end{claim}
    
    \begin{proofclaim}
    Figure~\ref{fig:claim_3/2} illustrates the different steps of the proof of the claim. The main steps are first to remove the colours of a set $S'$ in $H_i$ which then allows to remove the colours of a set $S''$ for vertices $v_{i+1}$ to $v_{i+\left \lceil \frac{k}{3} \right \rceil}$ and finally reach a colouring where no colour of $S''$ appears in $H_{i+\left \lceil \frac{k}{3} \right \rceil}$ (see the definitions of $S'$ and $S''$ below).

    \begin{figure}
        \begin{center}	
              \begin{tikzpicture}[thick,scale=0.9, every node/.style={transform shape}]
                \tikzset{vertex/.style = {circle,fill=black,minimum size=6pt, inner sep=0pt}}
                \tikzset{littlevertex/.style = {circle,fill=black,minimum size=4pt, inner sep=0pt}}
                \tikzset{edge/.style = {->,> = latex'}}

                \draw[pattern=north west lines] (0,0) rectangle ++(3,-0.5);
                \node[] (Gammai) at (-1.5,-0.25) {$\gamma_i$ : };
                \node[] (Hi) at (1.5,0.3) {$H_i$};
                \node[littlevertex, label={[label distance=4]-90:$v_1$}] (v1a) at (0.2,-0.25) {};
                \node[littlevertex, label={[label distance=4]-90:$v_i$}] (via) at (2.8,-0.25) {};
                
                \draw[] (3.3,0) rectangle ++(3.3,-0.5);
                \node[littlevertex, label={[label distance=4]-90:$v_{i+1}$}] (vip1a) at (3.5,-0.25) {};
                \node[] (dots1a) at (4.2,-0.25) {$\cdots$};
                \node[littlevertex, label={[label distance=4]-90:$v_j$}] (vja) at (4.9,-0.25) {};
                \node[] (dots2a) at (5.6,-0.25) {$\cdots$};
                \node[littlevertex, label={[label distance=4]-90:$v_{i+\left \lceil \frac{k}{3} \right \rceil}$}] (vik3a) at (6.3,-0.25) {};
(8.9,-2.25)
                \draw[] (6.9,0) rectangle ++(3,-0.5);
                \node[littlevertex, label={[label distance=4]-90:$v_n$}] (vna) at (9.7,-0.25) {};
                
                \draw[fill=white,white] (8.4,-0.25) ellipse (1 and 0.5);
                \draw[] (8.4,-0.25) ellipse (1 and 0.5);
                \node[] (cj) at (8.4,-0.25) {$ \neq c_j$};
                \node[] (Xj) at (8.4,0.65) {$X_j$};
                
                \draw[edge] (vja) to[in=135, out=35] (7.63, 0.1);
                \draw[edge] (vja) to[in=150, out=35] (8.4, 0.25);
                \draw[edge] (7.63, -0.6) to[in=-45, out=-135] (vja);

                \draw[pattern={Dots}] (0,-2+0.3) rectangle ++(3,-0.5);
                \node[] (Gammaip) at (-1.5,-2.25+0.3) {$\gamma_i'$ : };
                \node[littlevertex, label={[label distance=4]-90:$v_1$}] (v1b) at (0.2,-2.25+0.3) {};
                \node[littlevertex, label={[label distance=4]-90:$v_i$}] (vib) at (2.8,-2.25+0.3) {};
                
                \draw[] (3.3,-2+0.3) rectangle ++(3.3,-0.5);
                \node[littlevertex, label={[label distance=4]-90:$v_{i+1}$}] (vip1b) at (3.5,-2.25+0.3) {};
                \node[] (dots1b) at (4.2,-2.25+0.3) {$\cdots$};
                \node[littlevertex, label={[label distance=4]-90:$v_j$}] (vjb) at (4.9,-2.25+0.3) {};
                \node[] (dots2b) at (5.6,-2.25+0.3) {$\cdots$};
                \node[littlevertex, label={[label distance=4]-90:$v_{i+\left \lceil \frac{k}{3} \right \rceil}$}] (vik3b) at (6.3,-2.25+0.3) {};
(8.9,-2.25)
                \draw[] (6.9,-2+0.3) rectangle ++(3,-0.5);
                \node[littlevertex] (vippb) at (7.1,-2.25+0.3) {};
                \node[] (dots3b) at (8.4,-2.25+0.3) {$\cdots$};
                \node[littlevertex, label={[label distance=4]-90:$v_n$}] (vnb) at (9.7,-2.25+0.3) {};

                \draw[pattern={Dots}] (0,-4+0.7) rectangle ++(3,-0.5);
                \node[] (Etai) at (-1.5,-4.25+0.7) {$\eta_i$ : };
                \node[littlevertex, label={[label distance=4]-90:$v_1$}] (v1c) at (0.2,-4.25+0.7) {};
                \node[littlevertex, label={[label distance=4]-90:$v_i$}] (vic) at (2.8,-4.25+0.7) {};
                
                \draw[fill=lightgrey] (3.3,-4+0.7) rectangle ++(3.3,-0.5);
                \node[littlevertex, label={[label distance=4]-90:$v_{i+1}$}] (vip1c) at (3.5,-4.25+0.7) {};
                \node[] (dots1c) at (4.2,-4.25+0.7) {$\cdots$};
                \node[littlevertex, label={[label distance=4]-90:$v_j$}] (vjc) at (4.9,-4.25+0.7) {};
                \node[] (dots2c) at (5.6,-4.25+0.7) {$\cdots$};
                \node[littlevertex, label={[label distance=4]-90:$v_{i+\left \lceil \frac{k}{3} \right \rceil}$}] (vik3c) at (6.3,-4.25+0.7) {};
(8.9,-2.25)
                \draw[] (6.9,-4+0.7) rectangle ++(3,-0.5);
                \node[littlevertex] (vippc) at (7.1,-4.25+0.7) {};
                \node[] (dots3c) at (8.4,-4.25+0.7) {$\cdots$};
                \node[littlevertex, label={[label distance=4]-90:$v_n$}] (vnc) at (9.7,-4.25+0.7) {};

                \draw[fill=lightgrey] (0,-6+1.1) rectangle ++(3,-0.5);
                \node[] (Gammaik3) at (-1.5,-6.35+1.1) {$\gamma_{i+\left \lceil \frac{k}{3} \right \rceil}$ : };
                \node[littlevertex, label={[label distance=4]-90:$v_1$}] (v1d) at (0.2,-6.25+1.1) {};
                \node[littlevertex, label={[label distance=4]-90:$v_i$}] (vid) at (2.8,-6.25+1.1) {};
                
                \draw[fill=lightgrey] (3.3,-6+1.1) rectangle ++(3.3,-0.5);
                \node[littlevertex, label={[label distance=4]-90:$v_{i+1}$}] (vip1d) at (3.5,-6.25+1.1) {};
                \node[] (dots1d) at (4.2,-6.25+1.1) {$\cdots$};
                \node[littlevertex, label={[label distance=4]-90:$v_j$}] (vjd) at (4.9,-6.25+1.1) {};
                \node[] (dots2d) at (5.6,-6.25+1.1) {$\cdots$};
                \node[littlevertex, label={[label distance=4]-90:$v_{i+\left \lceil \frac{k}{3} \right \rceil}$}] (vik3d) at (6.3,-6.25+1.1) {};
(8.9,-2.25)
                \draw[] (6.9,-6+1.1) rectangle ++(3,-0.5);
                \node[littlevertex] (vippd) at (7.1,-6.25+1.1) {};
                \node[] (dots3d) at (8.4,-6.25+1.1) {$\cdots$};
                \node[littlevertex, label={[label distance=4]-90:$v_n$}] (vnd) at (9.7,-6.25+1.1) {};

                \draw[pattern=north west lines] (-2+0.8,-8+1.1) rectangle ++(1,-0.5);
                \node[] (legend1) at (-0.1+0.8,-8.25+1.1) {avoids $S$};
                \draw[pattern={Dots}] (2+0.8,-8+1.1) rectangle ++(1,-0.5);
                \node[] (legend2) at (3.95+0.8,-8.25+1.1) {avoids $S'$};
                \draw[fill=lightgrey] (6+0.8,-8+1.1) rectangle ++(1,-0.5);
                \node[] (legend3) at (8+0.8,-8.25+1.1) {avoids $S''$};
                \draw[] (-1.7,-7.6+1.1) rectangle ++(11.5,-1.3);
                
                \draw[->] (Gammai) -- (Gammaip) node [midway,xshift=-2em] () {$\frac{4k+12}{3}n$};
                \draw[->] (Gammaip) -- (Etai) node [midway,xshift=-2em] () {$\left \lceil \frac{k}{3} \right \rceil$};
                \draw[->] (Etai) -- (Gammaik3) node [midway,xshift=-2em] () {$\frac{4k+12}{3}n$};
              \end{tikzpicture}
          \caption{The redicolouring sequence between $\gamma_i$ and $\gamma_{i+\left \lceil \frac{k}{3} \right \rceil}$.}
          \label{fig:claim_3/2}
        \end{center}
    \end{figure}
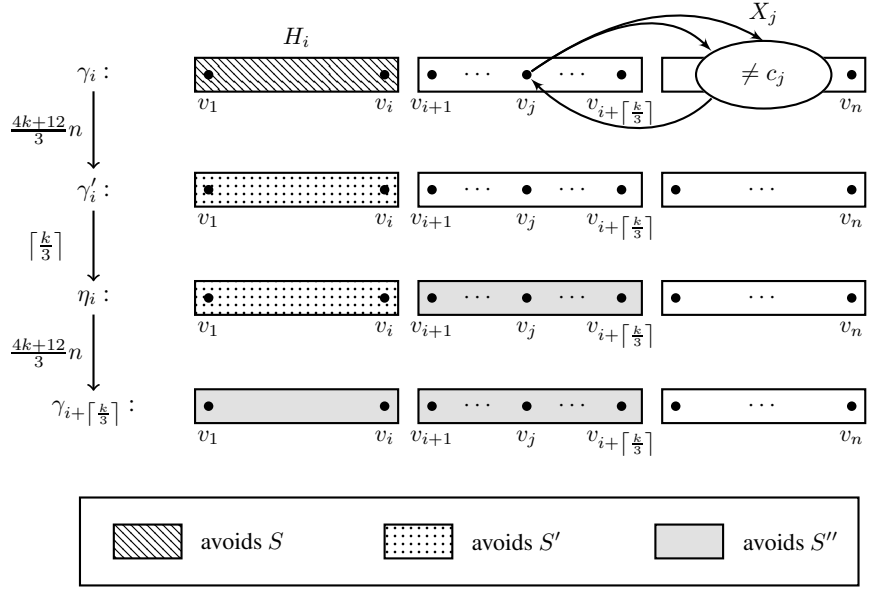
    
    Let $S$ be a set of colours of size exactly $\left \lceil \frac{k}{3} \right \rceil$ avoided by $\gamma_i$ on $H_i$.
    For each vertex $v_j$ in $\{v_{i+1},\dots,v_{i+\left \lceil \frac{k}{3} \right \rceil} \}$, we choose a colour $c_j$ so that each of the following holds:
    \begin{itemize}
        \item $c_j$ belongs to $L_j(v_j)$, and
        \item for each $\ell\in \{i+1,\dots,j-1\}$, $c_\ell$ is different from $c_j$.
    \end{itemize}
    Note that this is possible because $L_j$ is $\left \lceil \frac{k}{3} \right \rceil$-feasible.
    Now let $S'$ be the set $\{c_{i+1},\dots, c_{i+\lceil\frac{k}{3}\rceil}\}$. Observe that $|S'| = \left \lceil \frac{k}{3} \right \rceil$. By Lemma~\ref{lemma:avoid-colours}, there is, in $H_i$, an $L_i$-recolouring sequence  of length at most $\frac{4k+12}{3}n$ from $\gamma_i$ to some $\gamma_i'$ that avoids $S'$. This recolouring sequence extends to a redicolouring sequence in $D$ by Remark~\ref{remark:extendstoD}.
    In the obtained dicolouring, since $\gamma_i'$ avoids $S'$ on $H_i$, we can recolour successively $v_j$ with $c_j$ for all $i+1 \leq j \leq i+\left \lceil \frac{k}{3} \right \rceil$ (starting from $v_{i+1}$ and moving forward to $v_{i+\left \lceil \frac{k}{3} \right \rceil}$). This does not create any monochromatic cycle by choice of $c_j$. Let $\eta_i$ be the resulting dicolouring of $D$. Now we define  a list assignment $\Tilde{L}_i$ of $H_i$ as follows: 
    
    $$\Tilde{L}_i(v_j) = [k] \setminus \{ \eta_i(v) \mid v\in N(v_j) \cap \{v_{i+1},\dots,v_n\})) \text{ for all }j\in[i].$$
    
    Using the same arguments as we did for $L_i$, we get that $\Tilde{L}_i$ is $\left \lceil \frac{k}{3} \right \rceil$-feasible for $H_i$. Note that $\eta_i$ is an $\Tilde{L}_i$-colouring of $H_i$ that avoids $S'$. Let $S''$ be any set of $\left \lceil \frac{k}{3} \right \rceil$ colours disjoint from $S'$. By Lemma~\ref{lemma:avoid-colours}, there is, in $H_i$, an $\Tilde{L}_i$-recolouring sequence of length at most $\frac{4k+12}{3}n$ from $\eta_i$ to some $\eta_i'$ that avoids $S''$. This recolouring sequence extends directly to a redicolouring sequence in $D$. Since $S'$ is disjoint from $S''$, the obtained dicolouring is an $L_{i+\left \lceil \frac{k}{3} \right \rceil}$-colouring of $H_{i+\left \lceil \frac{k}{3} \right \rceil}$ that avoids at least $\left \lceil \frac{k}{3} \right \rceil$ colours in $H_{i+\left \lceil \frac{k}{3} \right \rceil}$. Hence we get a redicolouring sequence from $\gamma_i$ to the desired $\gamma_{i+\left \lceil \frac{k}{3} \right \rceil}$, in at most $\frac{8k + 24}{3}n + \left \lceil \frac{k}{3} \right \rceil$ steps.
    \end{proofclaim}
    
     Note that $\gamma_{\left \lceil \frac{k}{3} \right \rceil}$ (a dicolouring satisfying the assumptions of Claim~\ref{claim+k/3} for $i=\left \lceil \frac{k}{3} \right \rceil$) can be reached from $\alpha$ in less than $n$ steps: for all $j\in [\left \lceil \frac{k}{3} \right \rceil]$, choose a colour $c_j$ so that each of the following holds:
    \begin{itemize}
        \item $c_j$ belongs to $L_j(v_j)$, and
        \item for each $\ell\in [j-1]$, $c_\ell$ is different from $c_j$.
    \end{itemize} 
    Now we can recolour successively $v_1,\dots,v_{\left \lceil \frac{k}{3} \right \rceil}$ (in this order) to their corresponding colour in $\{c_1,\dots,c_{\left \lceil \frac{k}{3} \right \rceil} \}$. Then applying Claim~\ref{claim+k/3} iteratively at most $\left \lfloor \frac{n}{\left \lceil \frac{k}{3} \right \rceil} \right \rfloor \leq \frac{3n}{k}$ times, we get that there is a redicolouring sequence of length at most $n + \frac{3n}{k}\left(\frac{8k+24}{3}n + \frac{k}{3}\right)$ from $\alpha$ to a dicolouring $\alpha'$ of $D$ that is also a proper colouring of $G$.
     Note that there exists a constant $C_1$, independent of $k$, such that
     $n + \frac{3n}{k}\left(\frac{8k+24}{3}n + \frac{k}{3}\right) \leq C_1n^2$.
     
 \medskip
 
 Let $\alpha$ and $\beta$ be two $k$-dicolourings of $D$.
 As proved above, there is a redicolouring sequence of length at most $C_1n^2$ from $\alpha$ (resp. $\beta$) to a dicolouring $\alpha'$ (resp. $\beta'$) of $D$ that is also a proper colouring of $G$.
Since  ${\cal C}_k(G)$ has diameter at most $C_0n^2$, there is 
a recolouring sequence of $G$ of length at most $C_0n^2$ from $\alpha'$ to $\beta'$, which is also a redicolouring sequence of $D$ (since every proper colouring of $G$ is a dicolouring of $D$).
The union of those three sequences yields a redicolouring sequence from $\alpha$ to $\beta$ of length at most $(2C_1+C_0) n^2$.
\end{proof}

\begin{proof}[Proof of Theorem~\ref{thm:results_cycle_degeneracy}\ref{thm:results_cycle_degeneracy_statement_3}]
    Let $\alpha$ and $\beta$ be two $k$-dicolourings. Let us show by induction on the number of vertices that there exists a redicolouring sequence from $\alpha$ to $\beta$ where every vertex is recoloured at most $\delta_{c}^*(D)+1$ times. 
    
    If $n = 1$ the result is obviously true. Let $D$ be a digraph on at least two vertices,  let $u$ be a vertex such that $d_c(u) \leq \delta_{c}^*(D)$ and let $D' = D - u$.
    We denote by $\alpha'$ and $\beta'$ the dicolourings of $D'$ induced, respectively, by $\alpha$ and $\beta$.
    By induction and since $\delta^*_{c}(D') \leq \delta^*_{c}(D)$, there exists a redicolouring sequence from $\alpha'$ to $\beta'$ such that each vertex is recoloured at most $\delta_{c}^*(D)+1$ times.
    Now we consider the same recolouring steps to recolour $D$, starting from $\alpha$.
    If for some step $i$, it is not possible to recolour $v_i$ to $c_i$, this must be because $u$ is currently coloured $c_i$ and recolouring $v_i$ to $c_i$ would create a monochromatic directed cycle.
    By definition of cycle-degree, there exists a transversal $X$ of the directed cycles containing $u$, with  $|X|\leq \delta^*_c(D)$ and $u \notin X$. Since $k\geq 2\delta^*_{c}(D)+2$, there are at least $\delta^*_{c}(D)+2$ colours that do not appear in $X$.
    We choose $c$ among these colours so that $c$ does not appear in the next $\delta^*_{c}(D)+1$ recolourings of $X$, and we recolour $u$ with $c$. 
    
    Since $|X| \leq \delta_{c}^*(D)$ and since each vertex in $D'$ is recoloured at most $\delta_{c}^*(D)+1$ times, the total number of recolourings in $X$ is at most $\delta_{c}^*(D)(\delta_{c}^*(D)+1)$ in the redicolouring sequence obtained by induction. Hence, in this new redicolouring sequence, $u$ is recoloured at most $\delta_{c}^*(D)$ times. We finally have to set $u$ to its colour in $\beta$.  Doing so $u$ is recoloured at most $\delta_{c}^*(D)+1$ times. This concludes the proof.
\end{proof}

\section{Bounds on the diameter of \texorpdfstring{$\mathcal{D}_k(D)$}{Dk(D)} when \texorpdfstring{$k\geq \delta^*_c(D)+2$}{k >= dc*(D)+2}}
\label{section:nd1_mad}

This section is devoted to the proofs of Theorems~\ref{thm:nd1} and~\ref{thm:mad}, which are based on the proofs of~\cite{feghaliJCT147}.

\medskip 

In all the section, $f,g : \mathbb{N}^2 \xrightarrow{} \mathbb{N}$ are the functions defined as $f(s,t) = (s+1)!(2t)^s$ and $g(s,t)=2sf(s,t)+2s+1$ respectively.
The following is straightforward.
\begin{proposition}
    For every $s,t\in \mathbb{N}$, $s\neq 0$, the following inequalities hold:
    \begin{align}
        f(s,t) &\geq \sum_{q=1}^t \big(2(s+1)f(s-1,q) \big). \label{eq:first}\\
        g(s,t) &\geq 2f(s,t) + 2 + g(s-1,t).\label{eq:second}\\
        g(s,t) &= O_s(t^s).\label{eq:third}
    \end{align}
\end{proposition}

We now prove the following main lemma\footnote{In his original proof, Feghali claims to obtain, in the statement corresponding to the third item of Lemma~\ref{lemma:subgraph_H} for undirected graphs, a multiplicative factor $(s+1)$ instead of $(s+1)!$ (in the function $f$). Since we are not able to understand how the smaller factor is obtained in the original proof, we state our result with the larger factor of which we are sure of the correctness, which anyway is hidden in the asymptotic notation.}.

\begin{lemma}
    \label{lemma:subgraph_H}
    Let $D=(V,A)$ be a digraph, $(V_1,\dots,V_t)$ be a partition of $V$, and $s\geq 0$, $k\geq s+2$ be two integers. Let $h\in [t]$ be such that, for every $p\leq h$ and every $u\in V_p$, there exists $X_u \subseteq \bigcup_{i=p+1}^t V_i$ such that $|X_u| \leq s$ and $X_u$ intersects every directed cycle containing $u$ in $D - \bigcup_{i=1}^{p-1}V_i$.
    
    Then, for every $k$-dicolouring $\alpha$ of $D$ and for any colour $c\in [s+2]$, there exists a redicolouring sequence between $\alpha$ and some $k$-dicolouring $\beta$ such that:
    \begin{itemize}
        \item for every $v\in  \bigcup_{i=1}^h V_i, \beta(v) \neq c$ and $\beta(v) \leq s+2$,
        \item no vertex of $\bigcup_{i=h+1}^t V_i$ is recoloured, and
        \item each vertex in $\bigcup_{i=1}^h V_i$ is recoloured at most $f(s,h)$ times.
    \end{itemize}
\end{lemma}

\begin{proof}
    We proceed by induction on $s$. Assume first that $s=0$ and let $H$ the subdigraph of $D$ induced by $ \bigcup_{i=1}^h V_i$. We claim that $D$ does not contain any directed cycle which intersects $V(H)$. Indeed, if $D$ contains such a directed cycle $C$, let $q\in[h]$ be the smallest index such that $V(C) \cap V_q \neq \emptyset$, and let $u \in V(C) \cap V_q$. Then $C$ is a directed cycle containing $u$ in $D - \bigcup_{i=1}^{q-1}V_i$, so $X_u$ must intersect $C$. This yields a contradiction because $X_u = \emptyset$ (since $|X_u|=0$).
    Thus, since no directed cycle of $D$ intersects $V(H)$, in $\alpha$ we can recolour each vertex of $H$ with the colour $c'\in [2]$ different from $c$. Since $f(0,h) = 1$, we get the result.

    \medskip
    
    Assume now that $s>0$. Let $\mathscr{C}$ be the set of colours greater than $s+2$ and let $W$ be the set of vertices with colour $c$ or any colour $c'\in \mathscr{C}$ in $\alpha$. Formally, $\mathscr{C} = [s+3,k]$ and $W = \{ v\in V(D) \mid \alpha(v) = c \vee \alpha(v) \in \mathscr{C}\}$. If $k=s+2$, $\mathscr{C}$ is empty and $W$ is the set of vertices coloured $c$.
    
    Let $q\in [h]$ be the smallest index such that $V_q \cap W \neq \emptyset$. If such an index does not exist, then  we take $\alpha = \beta$ and we are done.
    Let $U_{q-1}=\bigcup_{i=1}^{q-1}V_i$ (when $q=1$ we let $U_0$ be the empty set). For each colour $a\in [s+2]$ different from $c$, we define $W_a$ as follows:
    \[
        W_a = \{ u \in W\cap V_q \mid \forall v\in X_u, \alpha(v) \neq a \}
    \]
    
    Observe that every vertex in $W_a$ is recolourable to colour $a$ in $D- U_{q-1}$ since $X_u$ intersects every directed cycle containing $u$ in $D- U_{q-1}$. Note also that every vertex $u\in W$ belongs to some $W_a$ (maybe to several) because $|X_u| \leq s$. Whenever a vertex belongs to several sets $W_a$, we remove it from one, so at the end the collection $(W_a)_{a\in [s+2],a\neq c}$ is a partition of $W$. Figure~\ref{fig:lemma_subgraph_H} illustrates the structure of $D$ dicoloured with $\alpha$.

    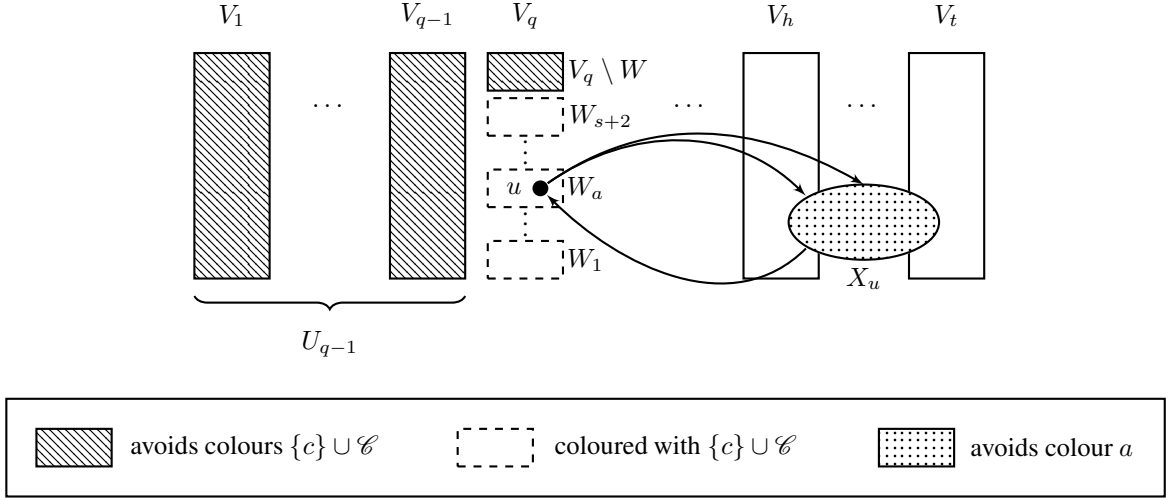
\begin{figure}
        \begin{center}	
              \begin{tikzpicture}[thick,scale=1, every node/.style={transform shape}]
                \tikzset{vertex/.style = {circle,fill=black,minimum size=6pt, inner sep=0pt}}
                \tikzset{littlevertex/.style = {circle,fill=black,minimum size=0pt, inner sep=0pt}}
                \tikzset{edge/.style = {->,> = latex'}}
            
                \draw[pattern=north west lines] (0,0) rectangle ++(1,-3);
                \node[] (V1) at (0.5,0.5) {$V_1$};
                \node[] (dots1) at (1.8,-0.7) {$\dots$};
                \draw[pattern=north west lines] (2.6,0) rectangle ++(1,-3);
                \node[] (Vqm1) at (3.1,0.5) {$V_{q-1}$};
                \draw[pattern=north west lines] (3.9,0) rectangle ++(1,-0.5);
                \node[] (Vq) at (4.4,0.5) {$V_q$};
                \node[] (VqmW) at (5.5,-0.25) {$V_q\setminus W$};
                
                \draw[dashed] (3.9,-0.6) rectangle ++(1,-0.5);
                \node[] (VqmW) at (5.4,-0.85) {$W_{s+2}$};
                
                \node[] (dots4) at (4.4,-1.22) {$\vdots$};
                
                \draw[dashed] (3.9,-1.55) rectangle ++(1,-0.5);
                \node[] (VqmW) at (5.2,-1.80) {$W_{a}$};
                \node[vertex, label=left:$u$] (u) at (4.6,-1.80) {};
                
                \node[] (dots4) at (4.4,-2.17) {$\vdots$};
                
                \draw[dashed] (3.9,-2.5) rectangle ++(1,-0.5);
                \node[] (VqmW) at (5.2,-2.75) {$W_{1}$};
                
                \node[] (dots2) at (6.6,-0.7) {$\dots$};
                
                \draw[] (7.3,0) rectangle ++(1,-3);
                \node[] (Vh) at (7.8,0.5) {$V_{h}$};
                
                \node[] (dots3) at (8.9,-0.7) {$\dots$};
                
                \draw[] (9.5,0) rectangle ++(1,-3);
                \node[] (Vh) at (10,0.5) {$V_{t}$};

                \draw[fill=white,white] (8.9,-2.25) ellipse (1 and 0.5);
                \draw[pattern={Dots}] (8.9,-2.25) ellipse (1 and 0.5);
                \node[] (Xu) at (8.9,-3) {$X_u$};
                
                \draw[decorate,decoration={brace,amplitude=5pt,mirror,raise=4ex}] (0,-2.6) -- (3.6,-2.6) node[midway,yshift=-3.6em]{$U_{q-1}$};

                \draw[edge] (u) to[in=135, out=35] (8.13, -1.9);
                \draw[edge] (u) to[in=150, out=35] (8.9, -1.75);
                \draw[edge] (8.13, -2.6) to[in=-45, out=-135] (u);
                
                \draw[pattern=north west lines] (-2.1,-5) rectangle ++(1,-0.5);
                \node[] (legend1) at (0.8,-5.25) {avoids colours $\{c\}\cup \mathscr{C}$};
                
                \draw[dashed] (3.5,-5) rectangle ++(1,-0.5);
                \node[] (legend2) at (6.4,-5.25) {coloured with $\{c\}\cup \mathscr{C}$};
                
                \draw[pattern={Dots}] (9.1,-5) rectangle ++(1,-0.5);
                \node[] (legend3) at (11.4,-5.25) {avoids colour $a$};

                \draw[] (-2.5,-4.6) rectangle ++(15.4,-1.3);
              
              \end{tikzpicture}
          \caption{The structure of the digraph $D$ dicoloured with $\alpha$, which we assume to be an ($s+2$)-dicolouring for clarity.  Note that $W_c$ does not exist and $X_u$ may intersect $V_{q+1}\cup \dots \cup V_{h-1}$.}
          \label{fig:lemma_subgraph_H}
        \end{center}
    \end{figure}
    
    \begin{claim}
        \label{claim:recolouring_Wa}
        Let $\phi$ be a $k$-dicolouring of $D$ such that:
        \begin{itemize}
            \item $\phi$ and $\alpha$ agree on $\bigcup_{i=q+1}^t V_i$ and $V_q \setminus W$, 
            \item $(\{c\} \cup \mathscr{C}) \cap \phi(U_{q-1}) = \emptyset$, and
            \item $\forall a\in [s+2], a\neq c$, either $\phi(W_a) = \{a\}$ or $\phi(W_a) \subseteq (\{c\} \cup \mathscr{C})$.
        \end{itemize}
        Then for every $a\in [s+2], a\neq c$ such that $\phi(W_a) \subseteq (\{c\} \cup \mathscr{C})$, there exists a redicolouring sequence between $\phi$ and a $k$-dicolouring $\psi$ such that:
        \begin{itemize}
            \item each vertex in $U_{q-1}$ is recoloured at most $2f(s-1,q-1)$ times,
            \item each vertex in $W_a$ is recoloured exactly once (to colour $a$),
            \item no vertex of $D - (U_{q-1} \cup W_a)$ is recoloured, and
            \item $(\{c\} \cup \mathscr{C}) \cap \psi(U_{q-1}) = \emptyset$.
        \end{itemize}
    \end{claim}
    \begin{proofclaim}
        By definition of $W_a$ and because $\phi$ and $\alpha$ agree on $\bigcup_{i=q+1}^t V_i$, note that every vertex $u\in W_a$ is recolourable to colour $a$ if $a \notin \phi(U_{q-1})$. Hence the key idea is to remove colour $a$ from $U_{q-1}$, then recolour every vertex in $W_a$ with $a$, and finally remove the colour $c$ from $U_{q-1}$ that we may have introduced. Along this process, we will never introduce any colour of $\mathscr{C}$. When $q=1$, note that we can just recolour every vertex in $W_a$ with $a$.

        \medskip

        Let $u_1,\dots,u_r$ be an ordering of $U_{q-1}$ such that the vertices in $V_p$ appear before the vertices in $V_{p'}$ for every $1\leq p' < p \leq q-1$. Whenever it is possible, in $\phi$, we recolour every vertex $u_1,\dots,u_r$ (in this order) with colour $c$. Let $\eta$ be the obtained dicolouring of $D$, and let $S$ be the set of vertices coloured $c$ in $\eta$. We define $D' = D - S$, $h' = q-1$ and $s' = s-1$. Also for every $p\in [t]$ we define $V_p' = V_p \setminus S$. Finally we define $\eta'$ as the induced dicolouring $\eta$ on $D'$. 

        \medskip
        
        Let us prove that, for every $p\leq h'$ and every $u\in V_p'$, the set of vertices $X_u' = X_u \setminus S$ satisfies $|X_u'|\leq |X_u|-1 \leq s'$ and intersects every directed cycle containing $u$ in $D' - \bigcup_{i=1}^{p-1}V_i'$. 
        
        First, since $u\in V_p'$, we know that $u$ has not been recoloured to $c$ in the previous process. It means that recolouring $u$ with $c$ creates a monochromatic directed cycle $C$. Moreover, since $c \notin \phi(U_{q-1})$ and by choice of the ordering $u_1,\dots,u_r$, we know that such a directed cycle $C$ is included in $\bigcup_{i=p}^tV_i$. By assumption on $X_u$, we have $X_u \cap V(C) \neq \emptyset$. Since $X_u \subseteq \bigcup_{i=p+1}^t V_i$ and $\phi(V(C) \setminus \{u\}) = \{c\}$ we deduce that $X_u \cap S \neq \emptyset$, which shows $|X_u'|\leq |X_u|-1 \leq s'$. 
        
        We now prove that  $X_u'$ intersects every directed cycle containing $u$ in $D' - \bigcup_{i=1}^{p-1}V_i'$. Let $C$ be such a directed cycle. Since $C$ is also a directed cycle in $D - \bigcup_{i=1}^{p-1}V_i$, we know that $C$ intersects $X_u$. We also know that $V(C) \cap S = \emptyset$ because $C$ is a directed cycle in $D'$. Hence, $C$ intersects $X_u \setminus S = X_u'$ as desired.

        \medskip

        By the remark above, we can apply the induction on Lemma~\ref{lemma:subgraph_H} with $D'$, $(V_1',\dots,V_t')$, $h'$, $s'$, $\eta'$ and $a$ playing the roles of  $D$, $(V_1,\dots,V_t)$, $h$, $s$, $\alpha$ and $c$ respectively. Hence, by induction, there exists a redicolouring sequence (which does not use colour $c$) in $D'$ from $\eta'$ to some dicolouring $\zeta'$ such that:
        \begin{itemize}
            \item for every $v\in \bigcup_{i=1}^{q-1}V_i', \zeta'(v) \notin (\{a\}\cup \mathscr{C})$,
            \item no vertex of $\bigcup_{i=q}^{t}V_i'$ is recoloured, and
            \item each vertex of $\bigcup_{i=1}^{q-1}V_i'$ is recoloured at most $f(s-1, q-1)$ times.
        \end{itemize}
        Since this redicolouring sequence does not use colour $c$, and because $\eta(S) = \{c\}$, it extends into a redicolouring sequence in $D$ between $\eta$ and $\zeta$, where $\zeta (u) = \zeta'(u)$ when  $u\in (U_{q-1} \setminus S)$, $\zeta(u) = c$ when $u\in S$ and $\zeta(u) = \phi(u)$ otherwise. Since $\zeta(v) \neq a$ for every vertex $v\in U_{q-1}$ and by choice of $W_a$, in $\zeta$ we can recolour every vertex in $W_a$ to colour $a$. Note that the vertices in $S\cap U_{q-1}$ have been recoloured exactly once (to colour $c$), which is less than $f(s-1,q-1)$.

        \medskip 
    
        We will now remove the colour $c$ we introduced in $U_{q-1}$. We use exactly the same process as before, swapping the roles of $c$ and $a$. So, whenever it is possible, starting with $u_1$ and moving forwards towards $u_r$, we recolour each vertex of $U_{q-1}$ with colour $a$. Let $\xi$ be the obtained dicolouring of $D$. We define $R= \{ v\in V(D) \mid \xi(v) = a \}$, $\Tilde{D} = D - R$, and $\Tilde{V}_i = V_i \setminus R$ for every $i\in [t]$. Finally let $\Tilde{\xi}$ be the induced dicolouring $\xi$ on $\Tilde{D}$. By induction, there exists a redicolouring sequence (which does not use colour $a$) in $\Tilde{D}$ from $\Tilde{\xi}$ to some dicolouring $\Tilde{\psi}$ such that:
        \begin{itemize}
            \item for every $v\in \bigcup_{i=1}^{q-1}\tilde{V}_i, \Tilde{\psi}(v) \neq c$,
            \item no vertex of $\bigcup_{i=q}^{t}\tilde{V}_i$ is recoloured, and
            \item each vertex of $\bigcup_{i=1}^{q-1}\tilde{V}_i$ is recoloured at most $f(s-1, q-1)$ times.
        \end{itemize}
        This gives, in $D$, a redicolouring sequence from $\xi$ to some dicolouring $\psi$ which does not use colour $c$ on $U_{q-1}$.
    
    \medskip

    Concatenating the redicolouring sequences we built, we conclude the existence of the desired redicolouring sequence from $\phi$ to $\psi$ in which every vertex in $U_{q-1}$ is recoloured at most $2 f(s-1,q-1)$ times, vertices in $W_a$ are recoloured exactly once to colour $a$, and the other vertices of $D$ are not recoloured.
    \end{proofclaim}

    \medskip

    Now we may apply Claim~\ref{claim:recolouring_Wa} on $\alpha$ (playing the role of $\phi$) to obtain a redicolouring sequence from $\alpha$ to a dicolouring $\alpha'$ (corresponding to $\psi$) in which $W_a$ has been recoloured to $a$ (for some fixed $a\in [s+2]$, $a\neq c$), and colours $\{c\}\cup \mathscr{C}$ do not appear in $\alpha'(U_{q-1})$.  
    Note that, in Claim~\ref{claim:recolouring_Wa}, the obtained dicolouring $\psi$ satisfies the assumptions on $\phi$. Thus, we may repeat this argument on $\alpha'$ to recolour $W_{a'}$ for some $a'\in [s+2], a'\notin \{c,a\}$. 
    
    Repeating this process for each colour $a\in [s+2]$ different from $c$, we obtain a redicolouring sequence between $\alpha$ and some dicolouring in which colours $\{c\}\cup \mathscr{C}$ do not appear in $U_{q-1} \cup V_q$ and such that each vertex of $U_{q-1}$ is recoloured at most $ (s+1) \cdot 2f(s-1,q-1) $ times and every vertex in $V_q$ is recoloured at most once.
    Since $f(s-1,q) \geq 1$ and $f(x,y)$ is non-decreasing in $y$, a fortiori every vertex of $\bigcup_{i=1}^qV_i$ is recoloured at most $2(s+1)f(s-1,q)$ times.
    
    \medskip

    We have shown above that, if colours $\{c\}\cup \mathscr{C}$ are not appearing in $\bigcup_{i=1}^{q-1} V_i$, then in at most $2(s+1)f(s-1,q)$ recolourings per vertex, we can also remove them from $V_{q}$ (and we do not recolour vertices in $\bigcup_{i=q+1}^{t} V_i$). Thus we can repeat this argument at most $h$ times to find a redicolouring sequence between $\alpha$ and a dicolouring $\beta$ in which colours $\{c\}\cup \mathscr{C}$ do not appear in $\bigcup_{i=1}^{h} V_i$. In this redicolouring sequence, the number of recolourings per vertex of $\bigcup_{i=1}^{h} V_i$ is at most
    
    \[ \sum_{q=1}^h \big( 2(s+1) f(s-1,q) \big) \leq f(s,h) \text{~~~~by Inequality~(\ref{eq:first}),}\]

    which concludes the proof.
\end{proof}

\begin{lemma}
    \label{lemma:diameter_mad}
    Let $D=(V,A)$ be a digraph on $n$ vertices and let $(V_1,\dots,V_t)$ be a partition of $V$ such that for every $p\in[t]$ and $u\in V_p$, there exists $X_u \subseteq \cup_{i=p+1}^t V_i$ such that $|X_u| \leq s$ and $X_u$ intersects every directed cycle containing $u$ in $D - \bigcup_{i=1}^{p-1}V_i$.
Then, for any $k\geq s+2$, $\mathcal{D}_{k}(D)$ has diameter at most $g(s,t)\cdot n$.
\end{lemma}
\begin{proof}
    We will show that, for any two $k$-dicolourings $\alpha,\beta$ of $D$, there exists a redicolouring sequence between them where each vertex is recoloured at most $g(s,t)$ times, showing the result. We proceed by induction on $s$. When $s=0$, $D$ is acyclic so we can directly recolour every vertex $v$ from $\alpha(v)$ to $\beta(v)$. Since $g(0,t) = 1$, we get the result.

    \medskip
    
    Assume now that $s>0$. By Lemma~\ref{lemma:subgraph_H}, there is a redicolouring sequence from $\alpha$ to an ($s+1$)-dicolouring $\Tilde{\alpha}$ in which each vertex of $D$ is recoloured at most $f(s,t)$ times (by taking $h=t$ and $c=s+2$). Symmetrically, we have a redicolouring sequence from $\beta$ to an ($s+1$)-dicolouring $\Tilde{\beta}$ in which each vertex of $D$ is recoloured at most $f(s,t)$ times. We will now find a redicolouring sequence between $\Tilde{\alpha}$ and $\Tilde{\beta}$.

    Let $v_1,\dots,v_n$ be an ordering of $V$ such that the vertices in $V_p$ appear before the vertices in $V_{p'}$ for every $1\leq p' < p \leq t$. In both $\Tilde{\alpha}$ and $\Tilde{\beta}$, starting with $v_1$ and moving forwards towards $v_n$, we recolour, whenever it is possible, each vertex of $V$ with colour $s+2$. This is done in at most two recolourings per vertex (one in both dicolourings). Let $\hat{\alpha}$ and $\hat{\beta}$ be the two obtained dicolourings. Observe that the vertices coloured $s+2$ in $\hat{\alpha}$ are exactly the vertices coloured $s+2$ in $\hat{\beta}$. We define $S = \{ v\in V \mid \hat{\alpha}(v) = s+2\}$ and $H = D - S$. Let $\hat{\alpha}_{|H}$ and $\hat{\beta}_{|H}$ be the dicolourings induced by $\hat{\alpha}$ and $\hat{\beta}$ on $H$, respectively. For each $p\in [t]$, let $V_p' = V_p\setminus S$. Observe that $(V_1',\dots,V_t')$ is a partition of $V(H)$ such that for every $p\in[t]$ and $u\in V_p'$, $X_u' = X_u \setminus S$ has size at most $s-1$ and intersects every directed cycle containing $u$ in $D - \bigcup_{i=1}^{p-1}V_i$ (the arguments are the same as in the proof of Lemma~\ref{lemma:subgraph_H}). Thus, by induction, there exists in $H$ a redicolouring sequence between $\hat{\alpha}_{|H}$ and $\hat{\beta}_{|H}$, using only colours in $[s+1]$ in which every vertex is recoloured at most $g(s-1,t)$ times. Since the vertices in $S$ are coloured $s+2$, this redicolouring sequence extends to $D$ and gives a redicolouring sequence between $\hat{\alpha}$ and $\hat{\beta}$. Thus, we have obtained a redicolouring sequence between $\alpha$ and $\beta$ in which the number of recolourings per vertex is at most $2f(s,t) + 2 + g(s-1,t) \leq g(s,t)$ by Inequality~(\ref{eq:second}).
    \end{proof}

We will now prove Theorems~\ref{thm:nd1} and~\ref{thm:mad} with Lemma~\ref{lemma:diameter_mad}.
\ndone*
\begin{proof}
    Take any cycle-degeneracy ordering $v_1,\dots,v_n$ of $D$, and set $V_i = \{v_i\}$ for every $i\in [n]$. Set $s=d$ and $t=n$, and the result follows directly from Lemma~\ref{lemma:diameter_mad} and Inequality~(\ref{eq:third}).
\end{proof}

\mad*
\begin{proof}
    Our goal is to find a partition $(V_1,\dots,V_{t(n)})$ of $V(D)$ such that $t(n) = O_{d,\epsilon}(\log n)$. Moreover, we need for every $p\in [t(n)]$ and every $u\in V_p$, that there exists $X_u \subseteq \bigcup_{i=p+1}^{t(n)}V_i$, $|X_u| \leq d-1$, that intersects every directed cycle containing $u$ in $D- \bigcup_{i=1}^{p-1}V_i$. If we find such a partition, then by Lemma~\ref{lemma:diameter_mad}, applied for $s=d-1$ and $t=t(n)$, we get that $\diam( \mathcal{D}_{k}(D) ) \leq g(d-1,t(n))\cdot n$, implying that $\diam( \mathcal{D}_{k}(D) ) = O_{d,\epsilon}(n (\log n)^{d-1})$ since $t(n) = O_{d,\epsilon}(\log n)$ and by Inequality~(\ref{eq:third}).

    \medskip 

    Let us guarantee the existence of such a partition.
    Let $H$ be any subdigraph of $D$ on $n_H$ vertices.
    For every vertex $u\in V(H)$, we let $X_u \subseteq V(H)\setminus \{u\}$ be a set of $d_c^H(u)$ vertices intersecting every directed cycle containing $u$ in $H$ (where $d_c^H(u)$ denotes the cycle-degree of $u$ in $H$). The existence of $X_u$ is guaranteed by definition of the cycle-degree.

    Then let $J=(V(H), F)$ be an auxiliary digraph, built from $H$, where $F = \{ uv \mid v \in X_u \}$. Let $S\subseteq V(H)$ be the set of all vertices $v$ with $d^+_J(v) \leq d-1$ (where $d_J^+(v)$ denotes the out-degree of $v$ in $J$). Then $|S| \geq \frac{\epsilon }{d}n_H$, for otherwise we have the following contradiction:
    \begin{align*}
        \Mad_c(D) \geq \Mad_c(H) &\geq \frac{1}{n_H} \sum_{v\in V(H)}d_c^H(v) = \frac{1}{n_H} \sum_{v\in V(H)}d_J^+(v)\\
        &\geq \frac{1}{n_H} \sum_{v\in V(H)\setminus S}d_J^+(v)\\
        &\geq \frac{1}{n_H} (n_H - |S|)d > (1 - \frac{\epsilon}{d})d = d- \epsilon,
    \end{align*}
    where in the last inequality we have used that $|S| < \frac{\epsilon n_H}{d}$.
    Now let us prove that $J[S]$ has an independent set $I$ of size at least $\frac{|S|}{2d-1}$. By choice of $S$, every subdigraph $J'$ of $J[S]$ satisfies $\Delta^+(J') \leq d-1$. Hence, for every such $J'$, we have $\sum_{v\in V(J')}(d_{J'}^+(v) + d_{J'}^-(v)) = 2|A(J')| \leq 2(d-1)|V(J')|$. In particular, this implies that $\UG(J[S])$ is ($2d-2$)-degenerate, and $\chi(\UG(J[S])) \leq 2d-1$. Take any proper ($2d-1$)-colouring of $\UG(J[S])$, its largest colour class is the desired $I$.

    Hence, we have shown that $H$ admits a set of vertices $I \subseteq V(H)$, of size at least $\frac{\epsilon}{(2d-1)d}n_H$, such that for every vertex $u\in I$ there exists $X_u \subseteq (V(H) \setminus I)$, $|X_u|\leq d-1$, that intersects every directed cycle of $H$ containing $u$.

    Since the remark above holds for every subdigraph $H$ of $D$, we can greedily construct  the desired partition $(V_1,\dots,V_{t(n)})$ by picking successively such a set $I$ in the digraph induced by the non-picked vertices. By construction, we get that $t$ satisfies the following recurrence:
    \[ t(i) \leq t\left(i-\frac{\epsilon i}{(2d-1)d} \right) + 1.\]

    Thus we have $t(n) \leq \log_{b}(n)$ where $b = \frac{1}{1-\frac{\epsilon}{(2d-1)d}}$, implying that 
    \[ t(n) \leq \frac{1}{-\log(1-\frac{\epsilon}{(2d-1)d})} \cdot \log(n) = O_{d,\epsilon}( \log(n)), \]
    which concludes the proof.
\end{proof}

\section{Bound on the diameter of \texorpdfstring{$\mathcal{D}_k(D)$}{Dk(D)} when \texorpdfstring{$k \geq \dic_g(D)+1$}{k >= Xg(D)+1}}
\label{section:digrundy_number}

This section is devoted to the proof of Theorem~\ref{thm:digrundy}.
\digrundy*
\begin{proof}
    Let $\alpha$ be any $k$-dicolouring of $D$ and $\beta$ be any $\dic(D)$-dicolouring of $D$. We will show by induction on $\dic(D)$ that there exists a redicolouring sequence of length at most $2\cdot \dic(D)\cdot n$ between $\alpha$ and $\beta$. The claimed result will then follow.
    If $\dic(D) = 1$, the result is clear since $D$ is acyclic.
    
    Starting from $\alpha$, whenever a vertex can be recoloured to colour $k$, we recolour it. Then we try to recolour the remaining vertices with colour $k-1$, and we repeat this process for every colour $k-1,\dots, 2$. At the end, the obtained dicolouring $\gamma$ is greedy (with colours ordered from $k$ to $1$). Actually, $\gamma$ is exactly the greedy dicolouring obtained from any ordering $v_1,\dots,v_n$ of $V(D)$ where $i<j$ whenever $\gamma(v_i) > \gamma(v_j)$. 

    Since $\gamma$ is a greedy dicolouring, and because $k\geq \dic_g(D) + 1$, colour 1 is not used in $\gamma$. This allows us to recolour every vertex of $V_1 = \{ v\in V(D) \mid \beta(v) = 1\}$ to colour 1 ($V_1\neq \emptyset$ since $\beta$ uses colours $[\dic(D)]$). If $\eta$ is the obtained dicolouring, then $\eta$ and $\beta$ agree on colour 1. Note also that, starting from $\alpha$, we reached $\eta$ by recolouring each vertex at most twice. Thus, the distance between $\alpha$ and $\eta$ in $\mathcal{D}_k(D)$ is at most $2n$. 
    
    Consider $H = D - V_1$. Since $\beta$ is an optimal dicolouring of $D$,  $\dic(H) = \dic(D) - 1$. Thus, by induction, there exists a redicolouring sequence between $\eta_{|H}$ and $\beta_{|H}$ (that is, the restrictions of $\eta$ and $\beta$, respectively, to $H$) of length at most $2(\dic(D) -1)n$, that does not use colour 1. This directly extends to a redicolouring sequence between $\eta$ and $\beta$ in $D$, which together with the redicolouring sequence between $\alpha$ and $\eta$ gives a redicolouring sequence between $\alpha$ and $\beta$ of length $2\cdot \dic(D) \cdot n$. 
\end{proof}

\section{Using the underlying graph to bound the diameter of \texorpdfstring{$\mathcal{D}_k(D)$}{Dk(D)}}
\label{section:general_bound_from_UG}

This section is devoted to the proof of Theorem~\ref{thm:general_bound_from_UG}.

\generalbound*

\begin{proof}
    Let $D=(V,A)$ be such a digraph, and let $\gamma$ be any proper $k$-colouring of $\UG(D)$. We will show, for any $k$-dicolouring $\alpha$ of $D$, that there is a redicolouring sequence between $\alpha$ and $\gamma$ of length at most $f(n)$, showing the result. Let $G_\alpha=(V,E)$ be the undirected graph where $E = \{ uv \mid uv\in A, \alpha(u) \neq \alpha(v) \}$. 

    First observe that $\alpha$ is a proper $k$-colouring of $G_\alpha$ by construction of $G_\alpha$. Note also that $\gamma$ is a proper $k$-colouring of $G_\alpha$ because $G_\alpha$ is a subgraph of $\UG(D)$. Moreover, since $G_\alpha$ is a subgraph of $\UG(D)$, by assumption on $\mathcal{G}$, we know that there exists a recolouring sequence between $\alpha$ and $\gamma$ in $G_\alpha$ of length at most $f(n)$.
    
    Next we show that every proper $k$-colouring of $G_\alpha$ is a dicolouring of $D$, implying that the recolouring sequence between $\alpha$ and $\gamma$ in $G_\alpha$ is also a redicolouring sequence between $\alpha$ and $\gamma$ in $D$. For purpose of contradiction, let us assume that $\beta$ is a proper $k$-colouring of $G_\alpha$ but $D$, coloured with $\beta$, contains a monochromatic directed cycle $C$. Then, by construction of $G_\alpha$, for each arc $xy$ of $C$, we must have $\alpha(x) = \alpha(y)$, for otherwise $xy$ would be a monochromatic edge in $G_\alpha$. This shows that $C$ is monochromatic in $D$ coloured with $\alpha$, a contradiction.
\end{proof}

Theorem~\ref{thm:general_bound_from_UG} directly extends to digraphs a number of known results about recolouring planar graphs. We discuss further these applications in Section~\ref{section:open_problems}.

\section{Case of digraphs with bounded \texorpdfstring{$\mathscr{D}$-width}{D-width}}
\label{section:directed_treewidth}

This section is devoted to the proof of Theorem~\ref{thm:dwidth}.

\dwidth*

The following claim can be easily deduced from the definition of a $\mathscr{D}$-decomposition.
\begin{claim}\label{claim:completeBag}
Let $(T,{\cal X}=(X_v)_{v \in V(T)})$ be a $\mathscr{D}$-decomposition of a digraph $D=(V,A)$ and $tt' \in E(T)$ such that $v \in X_{t'} \setminus X_t$. Then, $(T,{\cal X}'=(X'_v)_{v \in V(T)})$ such that $X'_u=X_u$ for all $u\neq t$ and $X'_t=X_t \cup \{v\}$ is a $\mathscr{D}$-decomposition of $D=(V,A)$. Moreover, if $|X_t|<|X_{t'}|$, $(T,{\cal X}')$ has the same width as  $(T,{\cal X})$.
\end{claim}

A $\mathscr{D}$-decomposition $(T,{\cal X})$ is {\it reduced} if, for every $tt' \in E(T)$, $X_t \setminus X_{t'}$ and $X_{t'} \setminus X_t$ are non-empty. It is easy to see that any digraph $D$ admits an optimal (i.e., of width $\dw(D)$) $\mathscr{D}$-decomposition which is reduced (indeed, if $X_t \subseteq X_{t'}$ for some edge  $tt' \in E(T)$, then contract this edge and remove $X_t$ from ${\cal X}$). 

A $\mathscr{D}$-decomposition $(T,{\cal X})$ of $\mathscr{D}$-width $k\geq 0$ is {\it full} if every bag has size exactly $k+1$. A $\mathscr{D}$-decomposition $(T,{\cal X})$ is {\it valid} if $|X_t \setminus X_{t'}|=|X_{t'}\setminus X_t|=1$ for every $tt' \in E(T)$. Note that any valid $\mathscr{D}$-decomposition is full and reduced. Note also that, if $(T,{\cal X})$ is valid and $t \in V(T)$ is a leaf of $T$, then there exists a (unique) vertex $v \in V$ that belongs only to the bag $X_t$. Such a vertex $v$ is called a {\it baby}.
\begin{lemma}\label{lem:valid_decomposition}
Every digraph $D=(V,A)$ admits a valid $\mathscr{D}$-decomposition of width $\dw(D)$.
\end{lemma}
\begin{proof}
Let $(T,{\cal X})$ be an optimal reduced $\mathscr{D}$-decomposition of $D=(V,A)$, which exists by the remark above the lemma. We will progressively modify $(T,{\cal X})$ in order to make it first full and then valid.

While the current decomposition is not full, let $tt' \in E(T)$ such that $|X_t| < |X_{t'}|=\dw(D)+1$ and let $v \in X_{t'} \setminus X_t$. Add $v$ to $X_t$. The obtained decomposition is still a $\mathscr{D}$-decomposition of width $\dw(D)$ by Claim~\ref{claim:completeBag}. Moreover, the updated decomposition remains reduced all along the process, as since $|X_t| < |X_{t'}|$ and the initial decomposition is reduced, $X_{t'}$ must contain another vertex $u \neq v$ with $u \notin X_t$. Eventually, the obtained decomposition $(T,{\cal X})$ becomes an optimal full $\mathscr{D}$-decomposition.

Now, while $(T,{\cal X})$ is not valid, let $tt' \in E(T)$, $x,y \in X_t \setminus X_{t'}$ and $u,v \in X_{t'} \setminus X_t$ (such an edge of $T$ and four distinct vertices of $V$ must exist since  $(T,{\cal X})$ is full and reduced but not valid). Then, add a new node $t''$ to $T$, with corresponding bag $X_{t''}=(X_{t'} \setminus \{u\}) \cup \{x\}$ and replace the edge $tt'$ in $T$ by the two edges $tt''$ and $t''t'$. Clearly, subdividing the edge $tt'$ by adding a bag $X_{t''}=X_{t'}$ still leads to an optimal full (but not reduced) $\mathscr{D}$-decomposition of the same width. Then, adding $x$ to $X_{t''}$ makes that $(T,{\cal X})$ remains a $\mathscr{D}$-decomposition (by the first statement of Claim~\ref{claim:completeBag}). Finally, we must prove that removing $u$ from $X_{t''}$ preserves the fact that we still have a $\mathscr{D}$-decomposition. Indeed, let $S$ be a strong subset whose support $T_S$  (before the subdivision) contains $tt'$ (clearly, the other strong subsets are not affected by the change in the decomposition). It must be because of some vertex in $z \in X_t \cap X_{t'}$ and so $z \in X_{t''}$. Therefore, the support $T_S$, obtained after the subdivision and the modifications to $X_{t''}$, contains both edges $tt''$ and $t''t'$, and therefore it remains connected. Note that, after the modifications, $(T,{\cal X})$ is still full and reduced.

Note that, after the application of each step as described above, either the maximum of $|X_t\setminus X_{t'}|$ over all edges $tt' \in E(T)$, or the number of edges $tt' \in E(T)$ that maximize $|X_t\setminus X_{t'}|$, strictly decreases, and none of these two quantities increases. Therefore, the process terminates, and eventually $(T,{\cal X})$ becomes an optimal valid $\mathscr{D}$-decomposition.
\end{proof}

Given a valid $\mathscr{D}$-decomposition $(T,{\cal X})$ of a digraph $D=(V,A)$, two vertices $u,v \in V$ are {\it parents}, denoted by $u \sim_p v$, if their supports $T_u$ and and $T_v$ (we use $T_v$ instead of $T_{\{v\}}$ for denoting the support of a single vertex $\{v\}$)  are vertex-disjoint and there is an edge $tt' \in E(T)$ with $t \in V(T_v)$ and $t' \in V(T_u)$. Let $\sim_{(T,{\cal X})}$ be the transitive closure of $\sim_p$.

\begin{lemma}\label{lem:classes}
Let $(T,{\cal X})$ be an optimal valid $\mathscr{D}$-decomposition of a digraph $D=(V,A)$. Then, $\sim_{(T,{\cal X})}$ defines an equivalence relation on $V$ which has exactly $\dw(D)+1$ classes. Moreover, the vertices of each class induce an acyclic subdigraph of $D$.
\end{lemma}
\begin{proof}
The facts that $\sim_{(T,{\cal X})}$ is well-defined and that there are $\dw(D)+1$ classes follow from the fact that $(T,{\cal X})$ is valid. Now, let $C$ be any equivalence class of $\sim_{(T,{\cal X})}$. For purpose of contradiction, let us assume that $D[C]$ contains a directed cycle $Q$. By definition of a $\mathscr{D}$-decomposition, the support $T_Q$ must induce a subtree of $T$. Since, by definition of $\sim_{(T,{\cal X})}$, the supports of the vertices of $Q$ are pairwise vertex-disjoint, the support $T_Q$ must consist precisely of the disjoint union of the supports of the vertices of $Q$, and hence $T_Q$ is not connected, a contradiction.
\end{proof}

Given a valid $\mathscr{D}$-decomposition $(T,{\cal X})$ of $D=(V,A)$, the corresponding equivalence relation $\sim_{(T,{\cal X})}$, and a subset $X \subseteq V$, a $k$-colouring $\alpha$ of $D$  (which is not necessarily a dicolouring) is {\it $X$-coherent with respect to $(T,{\cal X})$} if, for every $u,v \in X$ such that $u \sim_p v$, and for every $t \in V(T_v)$, $v$ is the unique vertex coloured $\alpha(v)$ in $X_t$, and $\alpha(u)=\alpha(v)$. In what follows, we will just say {\it $X$-coherent}, as the $\mathscr{D}$-decomposition will always be clear from the context.

\begin{claim}
Let $(T,{\cal X})$ be  a valid $\mathscr{D}$-decomposition of a digraph $D=(V,A)$ of width $k-1\geq 0$. Then, $D$ admits a $k$-colouring that is $V$-coherent.
\end{claim}
\begin{proofclaim}
    Consider any ordering $C_1,\dots,C_{k+1}$ of the classes of $\sim_{(T,{\cal X})}$. Let $\alpha$ be the colouring associating with each vertex $v$ the index $i$ such that $v\in C_i$.  Then $\alpha$ is $V$-coherent.
\end{proofclaim}

The following claim is straightforward, so we skip its proof.

\begin{claim}\label{claim:colourClasses}
Let $(T,{\cal X})$ be  a valid $\mathscr{D}$-decomposition of a digraph $D=(V,A)$ of width $k-1\geq 0$ and let $\alpha$ be a $k$-colouring of $D$ that is $V$-coherent. Then, the equivalence classes of $\sim_{(T,{\cal X})}$ are precisely $\alpha_1,\ldots,\alpha_k$ (the colours classes of $\alpha$). In particular, by Lemma~\ref{lem:classes}, $\alpha$ is a $k$-dicolouring.
\end{claim}

\begin{lemma}\label{lem:small-distance-coherent}
Let $D=(V,A)$ be an $n$-vertex digraph and let $(T,{\cal X})$ be  a valid $\mathscr{D}$-decomposition of $D=(V,A)$ of width $k-1=\dw(D)$. Let $\alpha$ and $\beta$ be two $(k+1)$-dicolourings of $D$ that are $V$-coherent. Then, $\alpha$ and $\beta$ are at distance at most $2n$ in $\mathcal{D}_{k+1}(D)$. 
\end{lemma}
\begin{proof}
We prove the existence of a redicolouring sequence from $\alpha$ to $\beta$ in $\mathcal{D}_{k+1}(D)$ that recolours each vertex at most twice.

Let $S_1,\dots,S_k$ be the equivalence classes of  $\sim_{(T,{\cal X})}$. By Claim~\ref{claim:colourClasses}, each $S_i$ corresponds exactly to one colour class of $\alpha$ and exactly one colour class of $\beta$. In particular, both $\alpha$ and $\beta$ use indeed $k$ colours (not necessarily the same).
Consider the undirected complete graph $H$ on $k$ vertices $x_1,\dots,x_k$, and the two colourings $\alpha_H,\beta_H$ of $H$ defined as $\{ \alpha_H(x_i) \} = \alpha(S_i)$ and $\{ \beta_H(x_i) \} = \beta(S_i)$. It is known (see~\cite[Lemma~5]{bonamyEJC69}) that there is a redicolouring sequence in $H$ between $\alpha_H$ and $\beta_H$ in which every vertex is recoloured at most twice. This directly extends to a redicolouring sequence between $\alpha$ and $\beta$ (when $x_i$ is recoloured with colour $c$ in $H$, we recolour every vertex in $S_i$ with $c$ in $D$). Note that this is indeed a redicolouring sequence because at each step of the sequence, in the corresponding colouring, every colour class is a subset of some $S_i$, which induces an acyclic subdigraph by Lemma~\ref{lem:classes}.
\end{proof}

Given a tree $T$ rooted in $r\in V(T)$ and two vertices $u,v$ of $T$, we say that $v$ is a \textit{descendant} of $u$ if $u$ belongs to the path between $r$ and $v$ in $T$.

\begin{lemma}\label{lem:subtree}
Let $D=(V,A)$ be an $n$-vertex digraph and $(T,{\cal X})$ be  a valid $\mathscr{D}$-decomposition of $D=(V,A)$ of width $k-1 \geq 0$. Let $T$ be rooted in $r \in V(T)$. Let $\alpha$ be a $(k+1)$-colouring of $D$ that is $(V\setminus X_r)$-coherent and let $c \in [k+1]$ such that $\alpha(v) \neq c$ for all $v \in X_r$. Then, for every $t,t'\in V(T)$ with $t'$ being a descendant of $t$, if there exists $v \in X_{t}$ with $\alpha(v)=c$, then $v$ is the unique vertex of $X_t$ coloured with $c$ and there exists a unique $u \in X_{t'}$ with $\alpha(u)=c$.  
\end{lemma}
\begin{proof}
For purpose of contradiction, let us assume that there exist $t,t' \in V(T)$ such that $t'$ is a descendant of $t$, a vertex in $X_{t}$ is coloured with $c$, and no vertex in $X_{t'}$ is coloured with $c$. Over all possible such pairs $\{t,t'\}$, we choose one such that the distance between $t$ and $t'$ in $T$ is minimum. Then $t'$ is a child of $t$. 

Let $v \in X_{t}$ such that $\alpha(v)=c$ and let $\{u\}= X_{t'} \setminus X_{t}$, where we have used that $(T,{\cal X})$ is valid. Note that $u\sim_{(T,{\cal X})} v$ and that, since $T_u$ is connected, $u \notin X_r$. Since $\alpha$ is $(V\setminus X_r)$-coherent, we must have $\alpha(u)=\alpha(v)=c$, a contradiction. 

The uniqueness of $u$ and $v$ comes from the fact that $u,v \notin X_r$ since $\alpha(v)=\alpha(u)=c$, and because by hypothesis $\alpha$ is $(V\setminus X_r)$-coherent.
\end{proof}

\begin{lemma}\label{lem:unifiy}
Let $D=(V,A)$ be an $n$-vertex digraph, let $(T,{\cal X})$ be  a valid $\mathscr{D}$-decomposition of $D=(V,A)$ of width $k-1$, and let $T$ be rooted in $r \in V(T)$. Let $y \in X_r$ and let $\alpha$ be a $(k+1)$-dicolouring of $D$ that is $(V\setminus (X_r\setminus \{y\}))$-coherent. Let $D'$ be the digraph obtained by identifying $y$ and all vertices $v$ such that $y \sim_{(T,{\cal X})} v$ and let $\alpha'$ be the dicolouring of $D$ arising from $\alpha$. If there exists a $(k+1)$-dicolouring $\beta'$ of $D'$ that can be reached from $\alpha'$ by recolouring each vertex at most once, then the $(k+1)$-dicolouring $\beta$ that naturally extends $\beta'$ to $D$ (i.e., $\beta(v)=\beta'(y)$ for all $v \sim_{(T,{\cal X})} y$) can be reached from $\alpha$ by recolouring each vertex at most once.
\end{lemma}
\begin{proof}
This follows from the fact that $\alpha$ is $(V\setminus (X_r\setminus \{y\}))$-coherent, and thus every bag $X_t$ has exactly one vertex of the colour of $y$, which is the colour of the unique vertex $v \sim_{(T,{\cal X})} y$ that belongs to $X_t$. Hence, to go from $\alpha$ to $\beta$, we follow the same redicolouring sequence for every vertex $u \nsim_{(T,{\cal X})} y$, and when we recolour $y$ in $D'$, we simply recolour every vertex $v \sim_{(T,{\cal X})} y$ with the same colour as $y$.
\end{proof}

\begin{claim}
\label{claim:properly_recoloured}
Let $D=(V,A)$ be an $n$-vertex digraph, let $(T,{\cal X})$ be  a valid $\mathscr{D}$-decomposition of $D=(V,A)$ of width $k-1$. Let $\alpha$ be a $(k+1)$-dicolouring of $D$ and $x\in V(D)$ be any vertex.
Let $c\in [k+1]$ be a colour such that, for every vertex $v\in \bigcup_{x\in X_t}X_t$, $\alpha(v) \neq c$. Then the ($k+1$)-colouring obtained from $\alpha$ by recolouring $x$ with $c$ is a dicolouring.
\end{claim}
\begin{proofclaim}
    Assume this is not the case, and recolouring $x$ with $c$ creates a monochromatic directed cycle $C$. Then since $D[V(C)]$ is strongly connected, the support $T_{V(C)}$ of $V(C)$ is a non-empty subtree. This implies that there exists $y\in V(C) \setminus \{x\}$ such that $y$ and $x$ belong to one same bag $X_t$. This contradicts the choice of $c$.
\end{proofclaim}

\begin{lemma}\label{lem:removeColour}
Let $D=(V,A)$ be an $n$-vertex digraph, let $(T,{\cal X})$ be  a valid $\mathscr{D}$-decomposition of $D=(V,A)$ of width $k-1$, and let $T$ be rooted in $r \in V(T)$. Let $\alpha$ be a $(k+1)$-dicolouring of $D$ that is $(V\setminus X_r)$-coherent and such that the colour $c$ does not appear in $X_r$,  i.e., there exists $c \in [k+1]$ such that $\alpha(v) \neq c$ for all $v \in X_r$. Then, there exists a $(k+1)$-dicolouring $\beta$ of $D$ and a redicolouring sequence $\alpha=\gamma_1,\dots,\gamma_\ell=\beta$ such that:
\begin{itemize}
    \item $\beta$ is $(V\setminus X_r)$-coherent,
    \item $\beta(v) \neq c$ for all $v \in V$,
    \item every vertex of $D\setminus X_r$ is recoloured at most once,
    \item no vertex of $X_r$ is recoloured, and
    \item if $x_i$ is the vertex recoloured between $\gamma_i$ and $\gamma_{i+1}$, then, for every vertex $v\in \bigcup_{x_i\in X_t}X_t$, $\gamma_i(v) \neq \gamma_{i+1}(x_i)$.
\end{itemize}
\end{lemma}
\begin{proof}
The proof is by induction on $k-1+|V(T)|$, where $k-1$ is the width of $(T,\mathcal{X})$. If $k-1= 0$, then $D$ is acyclic and colour $c$  can be eliminated by recolouring every vertex at most once with a same colour distinct from $c$. Note that the vertices in $X_r$ are not recoloured and the last condition holds trivially since every bag has size one.

If $|V(T)|=1$, the result holds trivially since the colour $c$ does not appear in $X_r$, so we may take $\beta = \alpha$. Hence, $r$ must have at least one child. Let us fix one child $v$ of $r$, and let $\{y\}=X_v \setminus X_r$. Let $T_v$ be the subtree of $T$ rooted in $v$ and let $D_v$ be the subdigraph of $D$ induced by $\bigcup_{t\in V(T_v)} X_t$.

\begin{itemize}
\item If $\alpha(y) \neq c$, then $c$ does not appear in $X_v$.  Let $(T_v,{\cal Y})=(T_v,\{X_t \mid t \in V(T_v)\})$ be the decomposition of $D_v$ obtained from $T$.  Let $D'_v$ be the digraph obtained from $D_v$ by identifying $y$ with all vertices of its class with respect to $(T_v,{\cal Y})$. Note that $(T_v,{\cal Y})$ is a full decomposition of $D'_v$. By contracting each edge $tt'\in E(T_v)$ such that $Y_t = Y_{t'}$ (in $D'$), we obtain $(T_v',{\cal Y}')$ a valid decomposition of $D'_v$. 

Note that $|V(T_v')|<|V(T)|$ and the width $(T_v',{\cal Y}')$ equals the width of $(T_v,{\cal X})$. Hence by induction there exists a $(k+1)$-dicolouring $\beta'_v$ of $D'_v$ that is $(V(D'_v)\setminus X_v)$-coherent and such that $\beta'_v(w) \neq c$ for all $w \in V(D'_v)$. Moreover, there is a redicolouring sequence from $\alpha'_v$, the dicolouring $\alpha$ restricted to $D'_v$, to $\beta'_v$ such that every vertex of $D'_v \setminus X_v$ is recoloured at most once, and vertices in $X_v$ are not recoloured. Note finally that, whenever a vertex $x$ is recoloured, it is recoloured with a colour that is not appearing in $ \bigcup_{x\in X_t}X_t$.

By Lemma~\ref{lem:unifiy},  there exists a $(k+1)$-dicolouring $\beta_v$ of $D_v$ that is $(V(D_v)\setminus X_v)$-coherent and such that $\beta_v(w) \neq c$ for every vertex $w \in V(D_v)$. Moreover, there is a redicolouring sequence $\gamma_v=(\gamma_1,\cdots,\gamma_{\ell})$ from $\gamma_1=\alpha_v$, the dicolouring $\alpha$ restricted to $D_v$, to $\gamma_{\ell}=\beta_v$ such that every vertex of $D_v \setminus X_v$ is recoloured at most once.  Furthermore, note that $\beta_v$ is indeed $(V(D_v) \setminus (X_v \setminus \{ y \})$-coherent because we identified $y$ with all vertices of its class in $D'_v$.

Along this redicolouring sequence $\gamma_v$, when a vertex $x$ is recoloured between $\gamma_i$ and $\gamma_{i+1}$, let us show that, for every vertex $z \in \bigcup_{x\in X_t}X_t$, $\gamma_i(z) \neq \gamma_{i+1}(x)$. Let us first assume that $x$ does not belong to the class of $y$ with respect to $(T_v,{\cal X})$. If $z \sim_{(T_v,{\cal X})} y$, then by induction, $\gamma_{i+1}(x) \neq \gamma_i(y)=\gamma_i(z)$. Otherwise ($z$ not in the class of $y$), by induction, we directly have that $\gamma_{i+1}(x) \neq \gamma_i(z)$. Second, let us assume that $x$ belongs to the class of $y$. Hence, $z \nsim_{(T_v,{\cal X})} y$ since $x$ and $z$ belong to a same bag. By induction, 

Finally, this redicolouring sequence in $D_v$ is indeed a redicolouring sequence in $D$ because of the property above and by Claim~\ref{claim:properly_recoloured}.

\item If $\alpha(y)=c$, by Lemma~\ref{lem:subtree} and because $(T,{\cal X})$ is $(V\setminus X_r)$-coherent, every bag in $T_v$ contains exactly one vertex coloured $c$ and the set of vertices coloured with $c$ is precisely $Y:=\{w \in X_t \mid t \in V(T_v),  w \sim_{(T,{\cal X})} y\}$. The reduced decomposition obtained from $(T_v,{\cal Y})=(T_v,\{X_t \setminus Y \mid t \in V(T_v)\})$ is a valid decomposition of $D'_v=D_v-Y$ of width $k-2$. Observe that it is full because we remove exactly one vertex from each bag. Let $c'$ be a colour that does not appear in $X_v \setminus \{y\}$. Observe that the width of $(T_v,{\cal Y})$ is at most $k-2$. Thus, by induction, there exists a $k$-dicolouring $\beta'$ of $D'_v$ that is $(V(D'_v)\setminus X_v)$-coherent and such that $\beta'(w) \notin \{c,c'\}$ for every vertex $w \in V(D'_v)$ and $\beta'$ can be obtained from $\alpha'_{v}$, the restriction of $\alpha$ to $D'_v$, by recolouring each vertex of $V(D'_v)\setminus X_v$ at most once. By then recolouring all vertices of $Y$ to $c'$, we obtain a $(k+1)$-dicolouring $\beta$ of $D_v$ that is $(V(D_v)\setminus (X_v \setminus \{y\}))$-coherent and such that $\beta(w) \neq c$ for every vertex $w \in V(D_v)$ and $\beta$ can be obtained from $\alpha_{v}$, the restriction of $\alpha$ to $D_v$, by recolouring each vertex of $D_v\setminus (X_v\setminus \{y\})$ at most once, as we wanted to prove.

Along the obtained redicolouring sequence, by induction, when a vertex $x$ that is not in the class $y$ is recoloured, it is recoloured with a colour different from the colours of $ \bigcup_{x\in X_t}(X_t\setminus Y)$. Since $Y$ is coloured with $c$, and no vertex is recoloured with $c$, it is recoloured with a colour different from the colours of $ \bigcup_{x\in X_t}(X_t\setminus Y)$. Finally, when $y' \sim_{(T,{\cal X})} y$ is recoloured, it is recoloured with $c'$ which is a colour that is not appearing in $D_v$. This shows the last property.

Finally, this redicolouring sequence in $D_v$ is indeed a redicolouring sequence in $D$ because of the property above and by Claim~\ref{claim:properly_recoloured}.
\end{itemize}

Repeating the process above for each child $v$ of $r$, we finally obtain a redicolouring sequence from $\alpha$ to some $(k+1)$-dicolouring $\beta$ such that $\beta(w) \neq c$ for all $w\in V$. Moreover, at each step, in both cases, no vertex of $X_r$ is recoloured, so no vertex of $X_r$ is recoloured all along the redicolouring sequence. Also, if $v$ and $v'$ are two children of $r$, and $x$ is a vertex of $D$ that belongs to $V(D_v) \cap V(D_{v'})$, then $x$ also belongs $X_r$, implying that $x$ is not recoloured. Thus, every vertex that is recoloured belongs to $V(D_v)$ for exactly one child $v$ of $r$, implying that, in the obtained redicolouring sequence, every vertex of $D\setminus X_r$ is recoloured at most once.

Finally, note that $\beta$ is $(V\setminus X_r)$-coherent because, in both cases, the obtained dicolouring $\beta_v$ is $(V \setminus \{ X_v \setminus \{y\} \})$-coherent.
\end{proof}

\begin{lemma}\label{lem:small-distance-general}
Let $D=(V,A)$ be an $n$-vertex digraph and $(T,{\cal X})$ be  a valid $\mathscr{D}$-decomposition of $D=(V,A)$ of width $k-1\geq \dw(D)$. For every $(k+1)$-dicolouring $\alpha$ of $D$ there exists a $V$-coherent $(k+1)$-colouring $\beta$ of $D$ such that $\alpha$ and $\beta$ are at distance at most $n^2$ in $\mathcal{D}_{k+1}(D)$. 
\end{lemma}
\begin{proof}
Let us root $T$ at $r \in V(T)$ arbitrarily. For a vertex $t \in V(T)$, let $T_t$ be the subtree of $T$ rooted at $t$ and let $D_t=D[\bigcup_{v \in V(T_t)} X_v]$. 

We will define inductively an ordering $(v_1,\ldots,v_n)$ on $V$ and a sequence $(\alpha=\gamma_0,\gamma_1,\ldots,\gamma_n=\beta)$ of $(k+1)$-dicolourings of $D$ such that $\gamma_i$ is $X_i$-coherent with $X_i=\{v_1,\ldots,v_i\}$ (set $X_0 = \emptyset$) for every $0 \leq i \leq n$ and such that it is possible to go from $\gamma_i$ to $\gamma_{i+1}$ by recolouring  every vertex of $X_{i}$ at most twice and $v_{i+1}$ at most once. Note that $\gamma_n=\beta$ is $V$-coherent.

First note that $\gamma_0$ is trivially $X_0$-coherent since $X_0 = \emptyset$.  

Let $i\geq 0$ and assume $(v_1,\ldots,v_i)$ and $(\gamma_1,\ldots,\gamma_i)$ that satisfy the above properties have already been defined. Let $v_{i+1} \in V \setminus X_i$ be any vertex that appears in some bag $X_t$  such that $V(D_t) \setminus X_t \subseteq X_i$ (if $i=0$ then $t$ is a leaf and $v_1$ is a baby). Note that ${\gamma_i}_{|D_t}$, the restriction of $\gamma_i$ to $D_t$, is a $(V(D_t) \setminus X_t)$-coherent dicolouring of $D_t$. Let $c$ be any colour that does not appear in $X_t$ coloured with $\gamma_i$. By Lemma~\ref{lem:removeColour}, there exists a $(k+1)$-dicolouring $\xi$ of $D_t$ that is $(V(D_t) \setminus X_t)$-coherent and such that $\xi(v) \neq c$ for every vertex $v \in V(D_t)$ and $\xi$ can be reached from ${\gamma_i}_{|D_t}$ by recolouring each vertex of $V(D_t) \setminus X_t \subseteq X_i$ at most once (when $i=0$, this sequence is empty). Note that the same recolouring sequence allows to go from $\gamma_i$ to the $(k+1)$-colouring $\gamma'_i$ (whose restriction to $D_t$ is~$\xi$) by recolouring each vertex of $X_i$ at most once and so that $\gamma'_i(v) \neq c$ for all $v \in V(D_t)$. Then, we can go from $\gamma'_i$ to $\gamma_{i+1}$ by recolouring each vertex $w \in V(D_t)$ such that $w \sim_{(T,{\cal X})} v_{i+1}$ with colour $c$ (note that all the vertices that are recoloured at this phase are in $X_{i}$ except $v_{i+1}$).  Then, $\gamma_{i+1}$ is $X_{i+1}$-coherent. Therefore, the induction properties hold for $i+1$.

At the end, we find the desired redicolouring sequence between $\alpha$ and $\beta$ in which the total number of recolourings  is at most:
\[ \sum\limits_{0 \leq i \leq n-1} (2|X_i|+1) = \sum\limits_{0 \leq i \leq n-1} (2i+1) = n^2,\]
which concludes the proof.
\end{proof}

Now the proof of Theorem~\ref{thm:dwidth} follows from Lemma~\ref{lem:small-distance-coherent} and Lemma~\ref{lem:small-distance-general}.

\dwidth*
\begin{proof}
    Take $\alpha$ and $\beta$ two $k$-dicolourings. Let $(T,\mathcal{X})$ be a valid $\mathscr{D}$-decomposition of $D$ of width $k-1 \geq \dw(D)$. By Lemma~\ref{lem:small-distance-coherent}, there is a redicolouring sequence from $\alpha$ (respectively $\beta$) to a dicolouring $\alpha'$ (respectively $\beta'$) that is $V$-coherent. Moreover, there is such a redicolouring sequence of length at most $n^2$. Then by Lemma~\ref{lem:small-distance-general}, there is a redicolouring sequence between $\alpha'$ and $\beta'$ of length at most $2n$. Altogether, this gives a redicolouring sequence between $\alpha$ and $\beta$ of length at most $2(n^2+n)$.
    
    Note that you can always build a valid $\mathscr{D}$-decomposition of width $k-1\geq \dw(D)$ unless $k\geq |V(D)|+1$, in which case the result is easy to prove.
\end{proof}

\section{Consequences on redicolouring planar digraphs  and open questions}
\label{section:open_problems}

In this work, we generalized several evidences for Cereceda's conjecture to digraphs. In particular, our results give more support to Conjecture~\ref{conj:directed_cereceda}. One can consider the following question which, if true, would imply Conjectures~\ref{conj:cereceda} and~\ref{conj:directed_cereceda}.
\begin{question}
    Let $k\in \mathbb{N}$ and $D$ be an $n$-vertex digraph such that $k \geq \delta^*_c(D)+2$. Is is true that $\diam(\mathcal{D}_k(D)) = O(n^2)$?
\end{question}
Using Proposition~\ref{prop:dtw_c_deg}, an analogue question is the following.
\begin{question}
    Let $k\in \mathbb{N}$ and $D$ be an $n$-vertex digraph such that $k \geq \dtw(D)+2$. Is is true that $\diam(\mathcal{D}_k(D)) = O(n^2)$?
\end{question}
The same question remains open when we replace  directed treewidth by  DAG-width or  Kelly-width. In every case, if true, it would give another generalisation of Theorem~\ref{thm:bonamy_bousquet_treewidth}.

Note that Conjecture~\ref{conj:directed_cereceda} implies Conjecture~\ref{conj:cereceda}. We ask if the converse is true.
\begin{question}
    Does proving Cereceda's conjecture for undirected  graphs imply its analogue in digraphs?
\end{question}

Finally we pose a few questions on redicolouring planar digraphs. The \emph{girth} of a graph is the size of a shortest cycle. The \emph{girth} of a digraph is the girth of its underlying graph. Using Theorem~\ref{thm:general_bound_from_UG}, we know that every result on recolouring  planar graphs extends to digraphs (up to a factor two).
In particular, a number of results from~\cite{bartierArxiv2021,bousquetJCTB155,dvorakEJC95,feghaliJCT147}, that we recap in Table~\ref{table:bounded_diameter_planar_graphs}, remain true on digraphs (the table is taken from~\cite{bartierArxiv2021}).

\renewcommand{\arraystretch}{1.5}
\begin{table}[H]
\centering
\begin{tabular}{ |c || c | c | c | c | c | c|} \hline
    \diagbox{$g$}{$k$} & 4 & 5 & 6 & 7 & 9 & $\geq 10$ \\ \hline \hline
    
    3 & $+\infty$ & $+\infty$ & $+\infty$ & $O(n^6)$~\cite{bousquetJCTB155} & $O(n^2)$~\cite{bousquetJCTB155} &$O(n)$~\cite{dvovrakEJC27} \\ \hline
    4 & $+\infty$ & $O(n^4)$~\cite{bousquetJCTB155} & $O(n\log^3 n)$~\cite{feghaliJCT147} & $O(n)$~\cite{dvorakEJC95} & - & - \\ \hline
    5 & $< +\infty$~\cite{bartierArxiv2021} & $O(n\log^2 n)$~\cite{feghaliJCT147} & - & - & - & - \\ \hline
    
    6 & $O(n^3)$~\cite{bousquetJCTB155} & $O(n)$~\cite{bartierArxiv2021} & - & - & - & - \\ \hline
    
    $\geq 7$ & $O(n\log n)$~\cite{feghaliJCT147} & - & - & - & - & - \\ \hline
    
\end{tabular}

\caption{Bound on the diameter of $\mathcal{C}_k(G)$ when $G$ is a planar graph with girth at least $g$. The value `$+\infty$' means that there exists a graph for which $\mathcal{C}_k(G)$ is not connected, and the value `$<+\infty$' means that $\mathcal{C}_k(G)$ is connected but no reasonable upper bound is known.}
\label{table:bounded_diameter_planar_graphs}
\end{table}

We ask if these results can be improved for oriented planar graphs (that means, planar digraphs with no digons). Using our results and some results of~\cite{bousquetArXiv23}, we obtain bounds on the diameter of $\mathcal{D}_k(\Vec{G})$ when $\Vec{G}$ is an oriented planar graph. We recap them in Table~\ref{table:bounded_diameter_oriented_planar_graphs}.
\begin{table}[H]
    \centering
    
    \begin{tabular}{ |c || c | c| c | c | c|} \hline
        \diagbox{$g$}{$k$} & 2 & 3 & 4 & 5 & 6  \\ \hline
        \hline
        
        3 & $+\infty$~\cite{bousquetArXiv23} & ? & $O(n^3)$ Th.~\ref{thm:nd1} & $O(n\log^3(n))$~Cor.~\ref{cor:mad} & $O(n)$~\cite{bousquetArXiv23} \\ \hline
        4 & ? & $O(n^2)$~\cite{bousquetArXiv23} & $O(n)$~\cite{bousquetArXiv23} & - & - \\ \hline
        5 & $< +\infty$~\cite{bousquetArXiv23} & $O(n\log (n))$ Cor.~\ref{cor:mad} & - & - & - \\ \hline
    \end{tabular}
    \caption{Bound on the diameter of $\mathcal{D}_k(\Vec{G})$ when $\Vec{G}$ is a oriented planar graph with girth at least $g$. The value `?' means that we do not know whether $\mathcal{D}_k(\Vec{G})$ is connected.}
    \label{table:bounded_diameter_oriented_planar_graphs}
\end{table}

We finally ask for the missing values in Table~\ref{table:bounded_diameter_oriented_planar_graphs}. 
\begin{question}\label{question:planar_3_mixing}
    Is every oriented planar graph 3-mixing?
\end{question}
\begin{question}
    Let $\Vec{G}$ be an oriented planar graph with $\girth(\Vec{G}) \geq 4$. Is it true that $\Vec{G}$ is 2-mixing? 
\end{question}
A famous conjecture by Neumann-Lara (see~\cite{bang2009}) states that every oriented planar graph has dichromatic number at most two. As a partial result, one may explore following question which, if true, implies that Question~\ref{question:planar_3_mixing} is a consequence of Neumann-Lara's conjecture.
\begin{question}
    Let $\Vec{G}$ be an oriented planar graph satisfying $\dic(G) \leq 2$. Is it true that $\Vec{G}$ is 3-mixing? 
\end{question}

\section*{Acknowledgments}

This research has been supported by the CAPES-Cofecub project Ma 1004/23, by the Inria Associated Team CANOE, by the research grant
    DIGRAPHS ANR-19-CE48-0013 and by the French government, through the EUR DS4H Investments in the Future project managed by the National Research Agency (ANR) with the reference number ANR-17-EURE-0004.

\bibliography{refs}

\end{document}